\documentclass[a4paper]{scrartcl}
\usepackage[]{amsmath,amssymb}
\usepackage{amsthm}
\usepackage{enumitem}
\usepackage[utf8]{inputenc}
\usepackage[T1]{fontenc}
\usepackage{lmodern}
\usepackage{fourier}
\usepackage{graphicx}
\usepackage{lmodern}
\usepackage{dsfont}
\usepackage{amsmath} % for theorem and proof environments
\usepackage{mathrsfs}
\usepackage{tikz}      % for graphics
\usepackage{tikz-cd}   % for commutative diagrams
\usepackage{mathtools}

\DeclareMathOperator{\oD}{d}
\def\D{\oD\!}
\def\dif{\,\D}

\usepackage{tikz}
\def\R{{\mathbf R}}

\def\Rps{\R_{+}^{*}}
\def\S{{\mathbf S}}

\def\N{{\mathbf N}}

\def\car#1{{\mathbf 1}}

\newcommand{\esp}[1]{{\mathbf E}\left[#1\right]}

\def\L{\mathbf L}
\def\P{\mathbb P}

\def\B{{\mathcal B}}%Tribu borélienne
\def\/{\,|\,} %Espérance conditionnelle
\def\car{\mathbf 1}
\usepackage{amsmath,amsfonts,amssymb,amsthm} % For math equations, theorems,
\newtheoremstyle{hermesexercises}{11pt}{11pt}{\normalfont}{0pt}{\scshape}{.--}{.5em}{}
{\theoremstyle{hermesexercises}
}
\newtheoremstyle{hermesremark}{11pt}{11pt}{\normalfont}{0pt}{\scshape}{.--}{.5em}{}
\newtheoremstyle{myhermesremark}{11pt}{11pt}{\normalfont}{0pt}{\scshape}{.--}{.5em}{}
{\theoremstyle{myhermesremark}
  \newtheorem{remark}{Remark}

}
\newtheoremstyle{hermestheorem}{11pt}{11pt}{\normalfont}{0pt}{\scshape}{.--}{.5em}{}
{\theoremstyle{hermestheorem}
\newtheorem{theoremT}{Theorem}
\newtheorem{definitionT}{Definition}
\newtheorem{corollaryT}[theoremT]{Corollary}
\newtheorem{lemmaT}[theoremT]{Lemma}}
\usepackage{tcolorbox}
\tcbuselibrary{skins,xparse,breakable}
\definecolor{fonce}{HTML}{323031}
\definecolor{cadet}{HTML}{084C61}
\definecolor{vertdeau}{HTML}{177E89}
\definecolor{ocre}{HTML}{FFC857}
\definecolor{creme}{HTML}{DB3A34}

\usepackage{bm}
\usepackage{leftindex}

\tcbset{textmarker/.style={%
skin=enhancedmiddle jigsaw,breakable,parbox=false, boxrule=0mm,leftrule=5mm,rightrule=5mm,boxsep=0mm,arc=0mm,outer arc=0mm, left=3mm,right=3mm,top=1mm,bottom=1mm,toptitle=1mm,bottomtitle=1mm,oversize}}

\tcbset{textmarker2/.style={%
skin=enhancedmiddle jigsaw,breakable,parbox=false, boxrule=0mm,leftrule=3mm,boxsep=0mm,arc=0mm,outer arc=0mm, left=3mm,right=3mm,top=1mm,bottom=1mm,toptitle=1mm,bottomtitle=1mm,oversize}}

\tcolorboxenvironment{lemmaT}{textmarker2,colback=ocre!5!white,colframe=ocre}
\newenvironment{lemma}{%
     \begin{lemmaT}}{\end{lemmaT}}

   \tcolorboxenvironment{definitionT}{textmarker2,colback=creme!5!white,colframe=creme}
   
      \tcolorboxenvironment{theoremT}{textmarker2,colback=vertdeau!5!white,colframe=vertdeau}
 \newenvironment{theorem}{%
     \begin{theoremT}}{\end{theoremT}}
      \tcolorboxenvironment{corollaryT}{textmarker2,colback=vertdeau!5!white,colframe=vertdeau}
 \newenvironment{proposition}{%
     \begin{theoremT}}{\end{theoremT}}
      \tcolorboxenvironment{corollaryT}{textmarker2,colback=vertdeau!5!white,colframe=vertdeau}
      \newenvironment{corollary}{%
     \begin{corollaryT}}{\end{corollaryT}}
\usepackage{enumitem}
\usetikzlibrary{%
  arrows,%
  calc,%
  fit,%
  patterns,%
  plotmarks,%
  shapes.geometric,%
  shapes.misc,%
  shapes.symbols,%
  shapes.arrows,%
  shapes.callouts,%
  shapes.multipart,%
  shapes.gates.logic.US,%
  shapes.gates.logic.IEC,%
  er,%
  automata,%
  backgrounds,%
  chains,%
  topaths,%
  trees,%
  petri,%
  mindmap,%
  matrix,%
  calendar,%
  folding,%
  fadings,%
  through,%
  positioning,%
  scopes,%
  decorations.fractals,%
  decorations.shapes,%
  decorations.text,%
  decorations.pathmorphing,%
  decorations.pathreplacing,%
  decorations.footprints,%
  decorations.markings,%
  shadows}
\newcommand{\iid}{\textnormal{i.i.d.\ }}
\newcommand{\lbra}{[\![}
\newcommand{\rbra}{]\!]}
\renewcommand{\esp}{\mathbb{E}}

\newcommand{\prob}{\mathbb{P}}
\newcommand{\ind}{\mathds{1}}

\newcommand{\lb}{\{}
\newcommand{\rb}{\}}

\newcommand\numberthis{\addtocounter{equation}{1}\tag{\theequation}}

\newcommand{\lipbis}{\mathrm{Lip}_{{[2]}}(\R)}
\def\eqdis{\overset{\mathrm{d}}{=}}

\def\Pi{\mathbf{P}^{1}}

\def\Pan{\mathbf{P}^{\alpha, \nu}}
\def\Lan{\mathscr{L}_{\alpha, \nu}}
\def\Dan{\mathbf{D}_{\alpha, \nu}}

\def\Pn{\mathbf{P}^{1,\nu}}

\def\cilipone{\mathcal{C}_{\text{Lip}_{1}}^{1}(\Eo)}
\def\cilog{\mathcal{C}_{\text{log}}^{1}(\Eos)}

\def\Epol{E_{\mathrm{pol}}}
\def\El{E_{\bm{\ell}}}
\def\Eo{E_{\bm{0}}}

\def\Eos{E^{*}_{\bm{0}}}
\def\m{\bm{m}}
\def\eqdis{\overset{\mathrm{d}}{=}}

\title{Stein's method for max-stable random vectors}
\author{B. Costacèque \and L. Decreusefond}
\date{2025}
\begin{document}
\maketitle{}
%\keywords{
\noindent{Keywords: Stein's method ; multivariate extreme value distributions ; Markov process ; generator approach}
%\msc{60G55}
\begin{abstract}
	Motivated by the omnipresence of extreme value distributions in 
limit theorems involving extremes of random processes, we adapt 
Stein's method to include these laws as possible target distributions.
We do so by using the generator approach of Stein's method, which is 
possible thanks to a recently introduced family of semi-groups. We study
the corresponding Stein solution and its properties when the working 
distance is either the smooth Wasserstein distance or the Kolmogorov distance.
We make use of those results to bound the distance between two max-stable 
random vectors, as well as to get a rate of convergence for the 
de Haan-LePage series in smooth Wasserstein distance.
\end{abstract}

%\begin{abstract}
%	We study the connections existing between max-infinitely divisible distributions and Poisson processes from the point of view of functional analysis. More precisely, we derive functional identities for the former by using well-known results of Poisson stochastic analysis. We also introduce a family of Markov semi-groups whose stationary measures are the multivariate max-stable distributions. As such, they constitute a tool to apply the generator approach of Stein's method to those distributions. We compute the generators of those semi-groups and deduce a Stein's identity characterizing max-stable distributions in higher dimensions. We solve the associated Stein's equation and study the properties of its solution. We then apply those results to obtain a rate of convergence in Wasserstein distance for the convergence of de Haan-LePage series. Additionally, we give a few functional identities associated to our semi-group, namely a Poincaré identity and commutation relations.   
%\end{abstract}

%\noindent{Keywords: Coupon collector's problem, Gumbel distribution, Stein's method, generator approach}

\noindent{Math Subject Classification: 60F05,60E07}
\setcounter{tocdepth}{1}
\tableofcontents

\section{Introduction}

The family of $\alpha$-stable distributions plays a central role in probability theory due to being the only possible non-degenerate limits for a renormalized sum of independent and identically distributed (i.i.d.) random variables. Since asymptotic results, like the central limit theorem and its generalization to non square-integrable distributions, are essential to applications in both statistics and probability, quantifying the speed of convergence to stable distributions is a recurring problem. Several techniques exist to deal with it, \textit{e.g.} Fourier analysis, the method of moments, \textit{etc.}, but they rarely work in a systematic way. Stein's method, introduced by C. Stein in \cite{Stein72} offers an alternative. It has successfully been used to obtain rates of convergence toward many types of target distributions, in particular the normal distribution and its higher dimensional counterparts (Gaussian vectors \cite{Götze91}, Gaussian processes and diffusions \cite{Barbour90}), the Poisson distribution (\cite{Chen74}) and the Poisson process (\cite{Barbour92, Chen04}), as well as $\alpha$-stable distributions (\cite{Chen22, Chen24, Coutin24}), among others.

However, few results of that kind exist when the limiting distribution is \textit{max-stable} or, equivalently, is an \textit{extreme value distribution} (EVD)). Those laws play a considerable role in extreme value theory, the subfield of probability theory focusing on estimating the frequency and intensity of extreme events, such as massive rainfalls, financial crashes, anomalous network traffic, \textit{etc.} A standard result about them is that they are equal to the distribution of the maximum of a certain Poisson process \cite{Resnick87}. For that reason, it is possible to apply the results developed in the aforementioned references \cite{Chen74, Barbour92, Chen04} to obtain rates of convergence toward EVD. For instance, if $X_{1},\dots, X_{n}$ are positive \iid random variables with tail function $F_{X}$, and whose tail function $\overline{F}_{X} = 1 - F_{X}$ is regularly varying with index $-\alpha$ for some $\alpha > 0$, \textit{i.e.} 
\[
\frac{\overline{F}_{X}(tx)}{\overline{F}_{X}(t)} = x^{-\alpha}, 
\]
then there exists a sequence $(a_{n})_{n \geq 0}$ of positive numbers such that $Z_{n} \coloneq a_{n}^{-1}\max(X_{1},\dots, X_{n})$ converges weakly to the Fréchet distribution $\mathcal{F}(\alpha)$. That result is actually a consequence of a more general statement \cite{Embrechts13}: the process of exceedances beyond $x > 0$
\[
\eta_{n,x} = \sum_{k=1}^{n}\delta_{a_{n}^{-1}X_{k}}\big([x, +\infty)\cap \cdot\big)
\]
converges weakly to a Poisson process $\eta_{x}$ with intensity measure $\alpha r^{-(\alpha+1)}$ on $[x, +\infty)$, where $\delta_{x}$ the Dirac measure in $x$. Assume one has a rate of convergence in total variation between those processes:
\[
\mathrm{d}_{\mathrm{TV}}(\eta_{n,x}, \eta_{x}) \leq u_{n,x},
\]
then by considering the event $\lb \eta[x, +\infty) = 0\rb$, one deduces immediately that
\[
\vert F_{Z_{n}}(x) - F_{Z}(x) \vert \leq u_{n,x},
\]
where $Z$ has the Fréchet distribution $\mathcal{F}(\alpha)$. However, one usually cannot deduce from this a rate of convergence in Kolmogorov distance for $Z_{n}$, because $u_{n,x}$ has no reason to be a bounded function of $x$. To bypass that limitation and obtain uniform rates of convergence, it seems reasonable to adapt Stein's method specifically to EVD.

Few instances of Stein's method being applied to EVD have been available so far. To the best of our knowledge, we can only cite \cite{Bartholome13}, which gives a uniform rate of convergence for the extreme value theorem in dimension 1. Even when using different techniques than Stein's method, that theorem has been quantified essentially in dimension 1 (\textit{e.g.} \cite{Hall79}, \cite{Cohen82}, \cite{Smith82}, \cite{deHaan96}), but much fewer results of that type exist in higher dimension (see \cite{Omey91} as well as A. Feidt's PhD dissertation \cite{Feidt13}). More generally, limit theorems involving multivariate max-stable distributions seem to have been seldom quantified. We propose to use Stein's method in order to provide a systematic framework to tackle those problems.

Our contributions are the following: we make use of the semi-groups $(\Pan_{t})_{t\geq 0}$ introduced in \cite{Costaceque24_mlti} to apply the generator approach of Stein's method to EVD, by first proving they provide a Stein solution and then by describing their main properties in the cases of the smooth Wasserstein and Kolmogorov distances. We use those results to bound the distance between two max-stable random vectors with different stability indexes and angular measures. We conclude by estimating the speed of convergence of the de Haan-LePage series in smooth Wasserstein distance.

To the best of our knowledge, those are the first results quantifying the speed of convergence of de Haan-LePage series. In the case of $\alpha$-stable distributions, the literature is richer: see for instance \cite{Ledoux96, Davydov99, Bentkus96, Bentkus00}, who managed to obtain a rate of convergence in total variation distance for all $\alpha \in (0,2)$. None relies on Stein's method to prove their result.

The rest of this paper is divided as follows: the first section recalls the basic notions and results we will rely of in the sequel, the second section extends $\Pan_{t}$ to test-functions not necessarily square-integrable with respect to the target distribution $\prob_{\alpha, \nu}$. It also introduces the associated Stein's solution and gives its main properties in the smooth Wasserstein and Kolmogorov cases. The third section is devoted to applications: bounding the distance between two max-stable random vectors and finding a rate of convergence for the de Haan-LePage series. 

\section{Notations and preliminaries}

The set of integers between $n$ and $m$ is denoted by $\lbra n,m \rbra$. Let $\bm{x} = (x^{1},\dots, x^{d})$ and $\bm{y} = (y^{1},\dots, y^{d})$ be two vectors in $\R^{d}$, with $x^{j} \leq y^{j}$ for all $j \in \lbra 1 ,d \rbra$. We set:
\[
[\bm{x}, \bm{y}] \coloneq \prod_{j=1}^{d}[x^{j}, y^{j}].
\]
Likewise, we take $[\bm{x}, \bm{y}) \coloneq \prod_{j=1}^{d}[x^{j}, y^{j})$. Let $\Eo$ be the set of vectors in $[\bm{0}, +\bm{\infty})$, minus the origin: 
\[
\Eo \coloneq [\bm{0}, +\bm{\infty}) \setminus \lb \bm{0} \rb.
\]
We will also work with the vectors $\bm{x}$ that are strictly greater than $\bm{\ell}$, in the sense that $x^{j} > 0$ for all $j \in \lbra 1,d \rbra$. We denote the set of such vectors by:
\[
\Eos \coloneq (\bm{0}, +\bm{\infty}).
\]
In the sequel, the notation $\bm{x} \leq \bm{y}$ means that the coordinates $x^{j}$ of $\bm{x}$ are less than or equal to their corresponding coordinates $y^{j}$ of $\bm{y}$, while $\bm{x} \nleq \bm{y}$ signifies that at least one coordinate of $\bm{x}$ is greater than its counterpart of $\bm{y}$. We denote the coordinate-wise maximum by:
\[
\bm{x}\oplus\bm{y} = \big(\max(x^{1}, y^{1}),\dots, \max(x^{d}, y^{d})\big)
\]
and
\[
\bm{x}\odot\bm{y} = \big(\min(x^{1}, y^{1}),\dots, \min(x^{d}, y^{d})\big).
\]
Besides $\max \bm{x} \coloneq \max(x^{1},\dots,x^{d})$ (respectively $\min \bm{x} \coloneq \min(x^{1},\dots,x^{d})$) denotes the greatest coordinate (respectively least) of $\bm{x}$. Consequently, it is always a scalar.

We say that a random vector $\bm{Z}$ is \textit{max-stable} if for all vectors $\bm{a}$, $\bm{b}$ in $\R_{+}^{d}$, there exists $\bm{c}, \bm{d} \in \R_{+}^{d}$, such that 
\begin{align}\label{prelim_ms_hyp}
\bm{a}\bm{Z}\oplus \bm{b}\bm{Z}' \eqdis \bm{c}\bm{Z} + \bm{d},
\end{align}
where $\bm{Z}'$ is an \iid copy of $\bm{Z}$. In \eqref{prelim_ms_hyp}, the sum and the multiplication between vectors are defined in a coordinate-wise way. A basic result in extreme value theory (see \cite{Resnick87} or \cite{Embrechts13} for instance) states that the marginals $Z^{j}$ of such a random vector $\bm{Z} = (Z^{1},\dots,Z^{d})$ are necessarily either Fréchet, Gumbel or Weibull random variables. The Fréchet distribution $\mathcal{F}(\alpha, \sigma)$ with shape parameter $\alpha > 0$ and scale parameter $\sigma > 0$ has c.d.f.
\begin{align}\label{prelim_Fréchet}
F(x) = \begin{cases}
e^{-\big(\frac{\sigma}{x}\big)^{\alpha}} &\quad \text{if\ } x \geq 0\\
0 &\quad \text{otherwise.}
\end{cases}
\end{align}
When $\sigma = 1$, we will simply note $\mathcal{F}(\alpha)$. In the sequel, we will assume that the $Z^{j}$ all have the same Fréchet distribution $\mathcal{F}(\alpha)$ for some $\alpha > 0$. When $\alpha = 1$, it is common to call such a random vector \textit{simple}. We will keep using this terminology for max-stable vectors whose marginals all have the same Fréchet $\mathcal{F}(\alpha)$ distribution. Simple max-stable random vectors have support on $\Eos$ and satisfy:
\begin{align}\label{prelim_stability}
\bm{a}\bm{Z}\oplus \bm{b}\bm{Z}' \eqdis \big(\bm{a}^{\alpha} + \bm{b}^{\alpha})\bm{Z},
\end{align}
where $\bm{x}^{\alpha}$ must be understood in a component-wise manner. We say a Radon measure $\mu$ on $\Eo$ possesses \textit{the $\alpha$-homogeneity property} if for all $t>0$:
\begin{align}\label{prelim_homo}
\mu\big(t^{\frac{1}{\alpha}}B\big) = t^{-1}\mu(B),\ B \in B(\Eo),
\end{align}
where $\B(\Eo)$ denotes the Borel $\sigma$-field of $\Eo$. Note that a Radon measure on $\Eo$ is $\sigma$-finite. We then have the most important theorem:

\begin{theorem}[de Haan-LePage representation]
Let $\alpha > 0 $ and $\bm{Z}$ a max-stable random vector with Fréchet $\mathcal{F}(\alpha)$ marginals. Then there exists $\eta = (\bm{y}_{i})_{i\geq 1}$ a Poisson process on $\Eo$ with a certain intensity measure $\mu$ such that the following equality in distribution holds: 
\begin{align}\label{prelim_LePage_1}
\bm{Z} \eqdis \bigoplus_{i = 1}^{\infty}\bm{y}_{i}.
\end{align}
\end{theorem}

In the sequel, $\mu$ will be called the \textit{exponent measure} of $\bm{Z}$. We refer to \cite{Last17}, \cite{Privault09} and the references therein for more about the Poisson process.

Thanks to the so-called \textit{polar decomposition}, it is possible to give more information about $\mu$. Fix a norm $\Vert \cdot \Vert$ on $\R^{d}$ (henceforth called the \textit{reference norm}) and set $\Epol \coloneq \Rps \times \S_{+}^{d-1}$, where $\S_{+}^{d-1}$ is the positive orthant of the sphere with respect to $\Vert\cdot\Vert$, \textit{i.e.}
\[
\S_{+}^{d-1} \coloneq \big\lb \bm{x} \in \R_{+}^{d},\ \Vert \bm{x} \Vert = 1\big\rb.
\]
For simplicity, we will assume that $\Vert\cdot\Vert$ is normalized so that $\S_{+}^{d-1} \subseteq [0,1]^{d}$. 
Define the transformation $T$
\begin{align*}
  \begin{array}{ccccc}
T & : & \Rps \times \S_{+}^{d-1} & \to & \Eo \\
& & (r, \bm{u}) & \mapsto & r\bm{u}^{\frac{1}{\alpha}}.
  \end{array}
\end{align*}
The polar decomposition of $\mu$, as stated in \cite{Resnick87} (proposition 5.11), implies that there exists $\nu$ a finite measure on $\S_{+}^{d-1}$ satisfying
\begin{equation}\label{prelim_moment_constraints}
\int_{\S^{d-1}_{+}}u^{j}\dif\nu(\bm{u}) = 1,\ j \in \lbra 1,d \rbra.
\end{equation}
and such that 
\begin{align}\label{prelim_polar}
\mu = T_{*}(\rho_{1}\otimes \nu)
\end{align}
where the right-hand side denotes the pushforward measure of $\rho_{1}\otimes \nu$ by $T$ and $\rho_{\alpha}$ is the measure on $\Rps$ defined by 
\begin{align}\label{prelim_rho}
\rho_{\alpha}[x, +\infty) \coloneq \frac{1}{x^{\alpha}}\cdotp
\end{align}
The previous result has the following consequence on the de Haan representation: there exists a marked Poisson process $\eta = ((r_{i}, \bm{u}_{i}))_{i\geq 1}$ on $\Epol$ such that
\begin{align}\label{prelim_LePage_2}
\bm{Z} \eqdis \bigoplus_{i=1}^{\infty}r_{i}\bm{u}_{i}^{\frac{1}{\alpha}}.
\end{align}
The scalar $\alpha$ is called the \textit{stability index} of $\bm{Z}$, while $\nu$ will be referred as the \textit{angular measure} of $\bm{Z}$. Since the distribution of a simple max-stable random vector is characterized equivalently by $\mu$ alone or $\alpha$ and $\nu$, we will parametrize it with either of them. We denote this by $\bm{Z} \sim \mathcal{MS}(\mu)$ and $\bm{Z} \sim \mathcal{MS}(\alpha, \nu)$ respectively.

\begin{equation}\label{prelim_moment_constraints}
\int_{\S^{d-1}_{+}}u^{j}\dif\nu(\bm{u}) = 1,\ j \in \lbra 1,d \rbra.
\end{equation}

In \cite{Costaceque24_mlti}, our set of test-functions was:
\[
\cilog = \big\lb f : \Eos \to \R,\ \exists C > 0,\ \vert f(\bm{x}) - f(\bm{y}) \vert \leq C\Vert \log \bm{x} - \log \bm{y} \Vert_{1} \forall \bm{x}, \bm{y} \in \Eos \big\rb,
\]
with $\Vert \bm{x} \Vert_{1} = \sum_{j=1}^{d}\vert x^{j} \vert$.
By contrast, our bounds will be expressed in terms of Kolmogorov distance $\mathrm{d}_{K}$ and Wasserstein distance $\mathrm{d}_{W}$, defined respectively as:
\[
\mathrm{d}_{K}(\bm{X}, \bm{Y}) \coloneq \underset{\bm{z} \in \Eos}{\sup}\vert F_{\bm{X}}(\bm{z}) - F_{\bm{Y}}(\bm{z})\vert = \Vert F_{\bm{X}} - F_{\bm{Y}} \Vert_{\infty},
\]
and
\[
\mathrm{d}_{W}(\bm{X}, \bm{Y}) \coloneq \underset{h \in \mathrm{Lip}_{1}(\R_{+}^{d})}{\sup}\vert \esp[h(\bm{X})] - \esp[h(\bm{Y})],
\]
where
\[
\mathrm{Lip}_{1}(\R_{+}^{d}) \coloneq \big\lb h : \R_{+}^{d} \to \R,\ \vert h(\bm{x}) - h(\bm{y}) \vert \leq \Vert \bm{x}-\bm{y} \Vert_{1}\ \text{for\ all\ } \bm{x}, \bm{y} \in \R^{d} \big\rb.
\]
For technical reasons, mainly due to the absence of density of max-stable distributions with respect to the Lebesgue measure in general, we will have to work with the smaller sets
\[
\cilipone \coloneq \mathcal{C}^{1}(\R_{+}^{d}) \cap \mathrm{Lip}_{1}(\R_{+}^{d}), 
\]
and 
\[
\mathrm{Lip}^{[2]}_{1}(\R_{+}^{d}) \coloneq \big\lb h : \R_{+}^{d} \to \R,\ h\ \text{and\ } \partial_{j}h\ \mathrm{are\ } 1\text{-Lipschitz}\ \forall j \in \lbra 1,d \rbra \big\rb.
\]
Notice that $\mathrm{Lip}^{[2]}_{1}(\R_{+}^{d})$ is a subset of $\cilipone$. We will refer to the associated metrics, especially the first one, as \textit{smooth Wasserstein distances}.

\section{Stein's equation and its solution}

Let $\bm{Z} \sim \mathcal{MS}(\alpha, \nu)$ be a max-stable random vector, with stability index $\alpha > 0$ and angular measure $\nu$. Assume that $\nu$ satisfies the moment constraint \eqref{prelim_moment_constraints}. Recall the definition of $\Pan_{t}$ when $\alpha$ is positive, as given in \cite{Costaceque24_mlti}:
\[
\Pan_{t}h(\bm{x}) = \esp\big[h\big(e^{-\frac{t}{\alpha}}\bm{x} \oplus (1-e^{-t})^{\frac{1}{\alpha}}\bm{Z}\big)\big],\ \bm{x} \in \R_{+}^{d},
\]
for $h \in \L^{p}(\P_{\bm{Z}})$ and $p \in [1,\infty]$. Its generator $\Lan$ is defined as:
$$
\Lan g \coloneq \underset{t \to 0^{+}}{\lim}\frac{\Pan_{t}g - g}{t}
$$
where the convergence is in norm $\Vert\cdot\Vert_{\L^{2}(\P_{\bm{Z}})}$. The operator $\Lan$ is equal to
$$
\Lan h(\bm{x}) = \frac{1}{\alpha}\langle \bm{x}, \nabla h(\bm{x}) \rangle + \Dan h(\bm{x})
$$
and $\Dan$ is defined as
$$
\Dan h(\bm{x}) = \int_{\El} \big(f(\bm{x}\oplus \bm{y}) - f(\bm{x})\big) \dif\mu(\bm{y}) = \frac{1}{\alpha}\int_{(\S_{+}^{d-1})^{1/\alpha}}\int_{\Rps}\big\langle r\bm{v}, \nabla f(\bm{x}\oplus r\bm{v}) \big\rangle_{\bm{x}} \frac{\alpha}{r^{\alpha + 1}}\dif r \dif\nu_{\alpha}(\bm{v}),
$$
$\nu_{\alpha}$ is the pushforward measure of $\nu$ by $\bm{x} \mapsto \bm{x}^{\alpha}$, $(\S_{+}^{d-1})^{1/\alpha}$ the set of elements of the form $\bm{v} = \bm{u}^{1/\alpha}$ for some $\bm{u} \in \S_{+}^{d-1}$ and
\[
\langle \bm{x}, \bm{y} \rangle_{\bm{z}} \coloneq \sum_{j=1}^{d}x^{j}y^{j}\ind_{\lb x^{j} \geq z^{j} \rb}.
\]

The semi-group $(\Pan_{t})_{t\geq 0}$ and its generator $\Lan$ are connected through the following relation.

\begin{lemma}
We have
\begin{align}\label{2_dif_Pan}
\frac{\dif}{\dif t}\Pan_{t}h(\bm{x}) = \Lan\Pan_{t}h(\bm{x})
\end{align}
for $h$ in $\cilipone$ and every $\bm{x} \in \R_{+}^{d}$, when $\alpha > 1$. 
%If $h$ is the indicator function of a segment $(-\bm{\infty}, \bm{z}]$ for some $\bm{z} \in \Eos$ and $\bm{x}$ belongs to
%\[
%A_{\bm{z}} = \big\lb \bm{x} \in \R_{+}^{d},\ \bm{x} < \bm{z}\ \mathrm{or\ } \bm{x} > \bm{z}\big\rb,
%\]
%then \eqref{2_dif_Pan} holds for every $\alpha > 0$.
\end{lemma}

\begin{proof}
Let $h \in \cilipone$. Adapting the proof given in \cite{Costaceque24_mlti} for the case $\alpha = 1$, we can prove for any $\alpha > 0$ that
\begin{align}\label{2_chaos_Pan}
\Pan_{t}h(\bm{x}) = e^{-\gamma_{t}\mu[\bm{0}, \bm{x}]^{c}}h\big(e^{-\frac{t}{\alpha}}\bm{x}\big) + \gamma_{t}e^{-\gamma_{t}\mu[\bm{0}, \bm{x}]^{c}}\int_{[\bm{0}, \bm{x}]^{c}}h\big(e^{-\frac{t}{\alpha}}(\bm{x}\oplus \bm{y})\big)\dif \mu(\bm{y}) + R_{t}(\bm{x}),
\end{align}
with $\bm{x} \in \Eos$, $\gamma_{t} \coloneq e^{t} - 1$ and
\[
R_{t}(\bm{x}) \coloneq \esp\Big[h\Big(e^{-\frac{t}{\alpha}}\big(\bm{x} \oplus \bigoplus_{i=1}^{N_{t, \bm{x}}}\bm{Y}_{i}\big)\Big)\ind_{\lb N_{t, \bm{x}} \geq 2 \rb}\Big]
\]
where $N_{t, \bm{x}} \sim \mathcal{P}(\gamma_{t}\mu[\bm{0}, \bm{x}]^{c})$ has the Poisson distribution with parameter $\gamma_{t}\mu[\bm{0}, \bm{x}]^{c}$, and the $\bm{Y}_{i}$ are $\iid$ random variables independent of $N_{t, \bm{x}}$ and whose distribution is given by 
\[
\prob(\bm{Y}_{1} \in A) = \frac{1}{\mu[\bm{0}, \bm{x}]^{c}}\mu(A),\ A \in \mathcal{B}([\bm{0}, \bm{x}]^{c}).
\] 
Differentiate equality \eqref{2_chaos_Pan} with respect to $t$ at $t=0$. The first term gives two parts of the generator:
\begin{align*}
\frac{\mathrm{d}}{\dif t}\Big\vert_{t=0}e^{-\gamma_{t}\mu[\bm{0},\bm{x}]^{c}}h(e^{-\frac{t}{\alpha}}\bm{x}) &= -\frac{1}{\alpha}\langle \bm{x},\nabla h(\bm{x}) \rangle - \mu[\bm{0},\bm{x}]^{c}h(\bm{x})\\
		&= -\frac{1}{\alpha}\langle \bm{x},\nabla h(\bm{x}) \rangle - \int_{[0,\bm{x}]^{c}}h(\bm{x})\dif \mu(\bm{y}).
\end{align*}
As for the first integral term, 
\begin{align*}
\frac{\mathrm{d}}{\dif t}\Big\vert_{t=0}\gamma_{t}e^{-\gamma_{t}\mu[\bm{0},\bm{x}]^{c}}\int_{[\bm{0}, \bm{x}]^{c}}h\big(e^{-\frac{t}{\alpha}}(\bm{x}\oplus & \bm{y})\big)\dif \mu(\bm{y})\\
	&= \underset{t \to 0}{\lim}\Big\lb e^{-\gamma_{t}\mu[\bm{0},\bm{x}]^{c}}\int_{[0,\bm{x}]^{c}}h\big(e^{-\frac{t}{\alpha}}(\bm{x}\oplus \bm{y})\big)\dif \mu(\bm{y}) \Big\rb\\
	&= \int_{[\bm{0},\bm{x}]^{c}}h(\bm{x}\oplus \bm{y})\dif \mu(\bm{y}).
\end{align*}
Notice that the integral term is differentiable with respect to $t$ because $\alpha$ is greater than 1. 
The remainder $R_{t}$ converges to $0$ at speed $o(t)$, so it does not contribute to the final result. Besides, as $\bm{y} \mapsto h(\bm{x}\oplus \bm{y}) - h(\bm{x})$ vanishes on $[\bm{0},\bm{x}]$, we find:
\[ 
\int_{[\bm{0},\bm{x}]^{c}}h(\bm{x}\oplus \bm{y})\dif \mu(\bm{y})-\int_{[\bm{0},\bm{x}]^{c}}h(\bm{x})\dif \mu(\bm{y}) = \int_{\Eo}\big(h(\bm{x}\oplus \bm{y}) - h(\bm{x})\big)\dif \mu(\bm{y}). 
\]
%Now let $h = \ind_{[-\bm{\infty}, \bm{z}]}$ and $\bm{x} \in A_{\bm{z}}$. Then $h$ is differentiable at $\bm{x}$ and its derivative is null. The same reasoning as above yields the desired result.
\end{proof}

The Markov semi-group $(\Pan_{t})_{t\geq 0}$ admits the law $\mathcal{MS}(\alpha, \nu)$ as its stationary distribution, and thus its generator $\Lan$ is a Stein operator for this distribution:
\begin{align}\label{2_LfZ}
\bm{Z} \sim \mathcal{MS}(\alpha, \nu) \Longrightarrow \esp[\Lan f(\bm{Z})] = 0,\ f \in \cilog.
\end{align}
We aim at using this generator to apply Stein's method to multivariate extreme value distribution. For every $h$ in the test-function set, this involves finding a solution $g_{h}$ to the Stein equation:
\begin{align}\label{2_Stein_eq}
\Lan g_{h}(\bm{x}) = h(\bm{x}) - \esp[h(\bm{Z})],\ \bm{x} \in \R_{+}^{d},
\end{align}
where $\bm{Z} \sim \mathcal{MS}(\alpha, \nu)$.

For a given $h$, the Stein solution $g_{h}$ of the equation $\Lan g_{h} = h - \esp[h(\bm{Z})]$ is formally given by 
\[
g_{h}(\bm{x}) = -\int_{0}^{\infty}\Pan_{t}h(\bm{x}).
\]
We check that definition makes sense whenever $h$ belongs to $\cilipone$ or the indicator function of $[\bm{0}, \bm{z}]$ for some $\bm{z} \in \Eos$.

\begin{proposition}
Let $\alpha > 0$ and $\bm{Z} \sim \mathcal{MS}(\alpha, \nu)$. Set $h^{*} = h - \esp[h(\bm{Z})]$ whenever $h$ is $\prob_{\alpha, \nu}$-integrable. We have
\[
\int_{0}^{\infty}\vert \Pan_{t}h(\bm{x}) - \esp[h(\bm{Z})] \vert\dif t < +\infty
\]
for every $\bm{x} \in \R_{+}^{d}$, and $h$ belongs to $\cilipone$ if $\alpha > 1$, or $h = \ind_{[\bm{0}, \bm{z}]}$ for some $\bm{z} \in \Eos$ and $\alpha > 0$. In the second case, the following relation holds:
\begin{align*}\label{2_psdo-inv_Kol}
-\int_{0}^{\infty} \Pan_{t}h^{*}(\bm{x}) \dif t = \alpha(\max \log \bm{x}\bm{z}^{-1}&)_{+}F_{\bm{Z}}(\bm{z})\\
		&- F_{\bm{Z}}(\bm{z})\int_{\alpha(\max \log \bm{x}\bm{z}^{-1})_{+}}^{\infty}\big(F_{\bm{Z}}(\bm{z})^{-e^{-t}} - 1\big)\dif t. \numberthis
\end{align*}
\end{proposition}

\begin{proof}
Assume $h$ belongs to $\cilipone$. Then
\begin{align*}
\vert \Pan_{t}h(\bm{x}) - \esp[h(\bm{Z})] \vert &\leq \esp\big[\Vert e^{-\frac{t}{\alpha}}\bm{x} \oplus (1-e^{-t})^{\frac{1}{\alpha}}\bm{Z} - \bm{Z}\Vert_{1}\big]\\
		&\leq \sum_{j=1}^{d}\esp\big[\vert e^{-\frac{t}{\alpha}}x^{j} - Z^{j}\vert\ind_{\lb Z^{j} \leq \gamma_{t}^{-1/\alpha}x^{j}\rb}\big] + d\big(1 - (1-e^{-t})^{\frac{1}{\alpha}}\big)\\
		&\leq \sum_{j=1}^{d}\esp\big[\vert e^{-\frac{t}{\alpha}}x^{j}- Z^{j}\vert^{\beta}\big]^{\frac{1}{\beta}}\prob(Z^{j} \leq \gamma_{t}^{-1/\alpha}x^{j})^{\frac{1}{\beta}} + d\big(1 - (1-e^{-t})^{\frac{1}{\alpha}}\big)\\
		&= \sum_{j=1}^{d}\esp\big[\vert e^{-\frac{t}{\alpha}}x^{j} - Z^{j}\vert^{\beta}\big]^{\frac{1}{\beta}}e^{-(\gamma_{t}(x^{j})^{-\alpha})^{\frac{1}{\beta}}} + d\big(1 - (1-e^{-t})^{\frac{1}{\alpha}}\big), 
\end{align*}
by Hölder's inequality, with $\beta = (\alpha + 1)/2$, so that $\beta \in (1,\alpha)$.
The first term is clearly integrable with respect to $t$, while the second one is equivalent to $e^{-t}/\alpha$ when $t$ goes to $\infty$. Thus $t\mapsto \Pan_{t}h(x) - \esp[h(Z)]$ is integrable w.r.t. Lebesgue measure.

If $h = h_{\bm{z}} = \ind_{[\bm{0}, \bm{z}]}$ for some $\bm{z} \in \Eos$, then $\Pan_{t}h_{\bm{z}}$ is well-defined because $\ind_{[\bm{0}, \bm{z}]}$ is bounded and measurable. The proof of identity \eqref{2_psdo-inv_Kol} goes as follows:
\begin{align*}
-\int_{0}^{\infty}\Pan_{t}h^{*}_{\bm{z}}(\bm{x})\dif t &= -\int_{0}^{\infty}\esp\big[\ind_{\lb e^{-t/\alpha}\bm{x} \oplus (1-e^{-t})^{1/\alpha} \leq \bm{z}\rb} - \prob(\bm{Z} \leq \bm{z})\big]\dif t\\
		&= -\int_{0}^{\infty} \prob((1-e^{-t})^{\frac{1}{\alpha}}\bm{Z}\leq \bm{z})\ind_{\lb e^{-t/\alpha}\bm{x} \leq \bm{z}\rb} - \prob(\bm{Z} \leq \bm{z})\dif t\\
		&= \alpha(\max \log \bm{x}\bm{z}^{-1})_{+}F_{\bm{Z}}(\bm{z}) - \int_{\alpha(\max \log \bm{x}\bm{z}^{-1})_{+}}^{\infty} F_{\bm{Z}}(\bm{z})^{1-e^{-t}} - F_{\bm{Z}}(\bm{z})\dif t.
\end{align*}
The last identity comes after noticing that
\begin{align*}
e^{-t/\alpha}\bm{x} \leq \bm{z} &\Longleftrightarrow \log \bm{x} \leq \frac{t}{\alpha}\bm{1} + \log \bm{z}\\
		&\Longleftrightarrow t \geq \alpha\log \frac{x^{j}}{z^{j}},\ \text{for\ all\ } j \in \lbra 1,d \rbra,
\end{align*}
which means that $t \geq \alpha (\max \log \bm{x}\bm{z}^{-1})_{+}$, the presence of the positive part coming from the fact that $t$ is always non-negative. The relation 
\[\
\prob\big((1-e^{-t})^{\frac{1}{\alpha}}\bm{Z} \leq \bm{z}\big) = e^{-\mu[\bm{0}, (1-e^{-t})^{-\frac{1}{\alpha}}\bm{z}]^{c}} = F_{\bm{Z}}(\bm{z})^{1-e^{-t}}
\]
stems from the homogeneity property of the exponent measure $\mu$ of $\bm{Z}$.
\end{proof}

We prove that $g_{h}$ is indeed the Stein solution of \eqref{2_Stein_eq}.

\begin{proposition}\label{2_Ln-1}

Let $\bm{Z} \sim	\mathcal{MS}(\alpha, \nu)$ and $h \in \cilipone$. Then $g_{h}$ is a solution of \eqref{2_Stein_eq}:
\begin{align}\label{2_Stein_sol_Wass}
\Lan g_{h}(\bm{x}) = h(\bm{x}) - \esp[h(\bm{Z})],
\end{align}
for all $\bm{x} \in \R_{+}^{d}$.
Furthermore, for any $\bm{z} \in \Eos$, the function $g_{\bm{z}}$ satisfies \eqref{2_Stein_eq} in the sense that:
\begin{align}\label{2_Stein_sol_Kol}
\Lan g_{\bm{z}}(\bm{x}) = \ind_{[\bm{0}, \bm{z}]}(\bm{x}) - F_{\bm{Z}}(\bm{z}),
\end{align}
for all $\bm{x} \in \R_{+}^{d}$ such that $x^{j} \neq z^{j}$ for every $j \in \lbra 1,d \rbra$.

\end{proposition}

\begin{proof}

Let $h \in \cilipone$. Thanks to the previous lemma, we have 
\[
\Lan\Pan_{t}g_{h}(\bm{x}) = \frac{\dif}{\dif t}\Pan_{t}g_{h}(\bm{x})
\]
Besides, 
\[
\Pan_{t}g_{h}(\bm{x}) = -\int_{0}^{\infty}\Pan_{t+s}h^{*}(\bm{x}) \dif s
\]
thanks to Fubini's theorem and the bound
\[
\big\vert h^{*}\big(e^{-\frac{t+s}{\alpha}}\bm{x} \oplus (1-e^{-(t+s)})^{\frac{1}{\alpha}}\bm{w}\big)\big\vert \leq \sum_{j=1}^{d}\vert e^{-\frac{t+s}{\alpha}}x^{j} - w^{j}\vert\ind_{\lb w^{j} \leq \gamma_{t+s}^{-1/\alpha}x^{j}\rb} + d\big(1 - (1-e^{-(t+s)})^{\frac{1}{\alpha}}\big),
\]
true $\dif s \otimes \dif\prob_{\bm{Z}}(\bm{w})$-a.s. Another dominated convergence argument yields that $\Pan_{t}h^{*}$ is differentiable with respect to each argument, with $j$-th partial derivative equal to
\[
\partial_{j}\Pan_{t}h^{*}(\bm{x}) = e^{-\frac{t}{\alpha}}\esp\big[(\partial_{j}h^{*})\big(e^{-\frac{t}{\alpha}}\bm{x} \oplus (1-e^{-t})^{\frac{1}{\alpha}}\bm{Z}\big)\ind_{\lb x^{j} \geq \gamma_{t}^{1/\alpha}Z^{j}\rb}\big]
\]
so that
\begin{align}\label{2_Pan_partl_bnd}
\vert \partial_{j}\Pan_{t+s}h^{*}(\bm{x}) \vert \leq e^{-\frac{t+s}{\alpha}}\esp\big[\ind_{\lb x^{j} \geq \gamma_{t+s}^{1/\alpha}Z^{j}\rb}\big] \leq e^{-\frac{s}{\alpha}}e^{-\gamma_{s}(x^{j})^{-\alpha}}.
\end{align}
On the other hand, we have 
\begin{align*}
\Dan \Pan_{t+s}h^{*}(\bm{x}) &= \int_{[\bm{0}, \bm{x}]^{c}} \big\langle \bm{y}, \nabla (\Pan_{t+s}h^{*})(\bm{x}\oplus \bm{y})\big\rangle_{\bm{x}} \dif \mu(\bm{y})\\
		&\leq \sum_{j=1}^{d} \int_{\lb y^{j} > x^{j}\rb} e^{-\frac{s}{\alpha}}e^{-\gamma_{s}(x^{j}\oplus y^{j})^{-\alpha}} \dif \mu(\bm{y}).
\end{align*}
When $x^{j} = 0$, the integral converges thanks to the exponential term in $y^{j}$. Otherwise, we can simply bound that term by $1$, since $\mu(\lb y^{j} > x^{j}\rb)$ is finite whenever $x^{j}$ is positive. Consequently, there exists a positive and $\dif s$-integrable function on $\R_{+}$ such that
\[
\vert \Lan \Pan_{t+s}h^{*}(\bm{x}) \vert \leq g(x^{j}, s),
\]
which implies that
\[
\Lan g_{h}(\bm{x}) = -\int_{0}^{\infty} \Lan\Pan_{s}h^{*}(\bm{x}) \dif s = -\int_{0}^{\infty} \frac{\dif}{\dif s}\Pan_{s}h(\bm{x}) \dif s = h^{*}(\bm{x})
\]
due to the ergodicity property of $(\Pan_{t})_{t\geq 0}$. 

In the case of the Kolmogorov distance, we use the two following facts: first the commutation rule
\[
\Dan \Pan_{t} = e^{-t}\Pan_{t}\Dan
\] 
and second, that $h_{\bm{z}}$ is an eigenvector for $\Dan$ with eigenvalue $-\mu[\bm{0},\bm{z}]^{c}$. Thus, by Fubini's theorem, we find
\[
\Dan g_{\bm{z}}(\bm{x}) = \mu[\bm{0},\bm{z}]^{c}\int_{0}^{\infty}e^{-t}\Pan_{t}h_{\bm{z}}(\bm{x}) \dif t.
\]
Consequently, since $F_{\bm{Z}}(\bm{z}) = e^{-\mu[\bm{0}, \bm{z}]}$,
\begin{align*}
\Dan g_{\bm{z}}(\bm{x}) &= \mu[\bm{0},\bm{z}]^{c}\int_{0}^{\infty}e^{-t}\prob\big((1-e^{-t})^{\frac{1}{\alpha}}\bm{Z} \leq \bm{z}\big)\ind_{\lb t \geq \alpha(\max \log \bm{x}\bm{z}^{-1})_{+}\rb} \dif t\\
		&= \mu[\bm{0},\bm{z}]^{c}\int_{\alpha(\max \log \bm{x}\bm{z}^{-1})_{+}}^{\infty}e^{-t}F_{\bm{Z}}(\bm{z})^{1-e^{-t}} \dif t\\
		&= \mu[\bm{0},\bm{z}]^{c}\int_{\alpha(\max \log \bm{x}\bm{z}^{-1})_{+}}^{\infty}e^{-t}e^{-\mu[\bm{0},\bm{z}]^{c}(1-e^{-t})} \dif t.
\end{align*}
Assume that the highest coordinate of $\log \bm{x}\bm{z}^{-1}$ is unique and has index $j$. Then we find:
\begin{align}\label{2_g_z_Dan}
\Dan g_{\bm{z}}(\bm{x}) = 
		\begin{cases}
		F_{\bm{Z}}(\bm{z})^{1 - \big(\frac{z^{j}}{x^{j}}\big)^{\alpha}} - F_{\bm{Z}}(\bm{z})  &\quad \text{if\ } x^{j} > z^{j}\\
		1 - F_{\bm{Z}}(\bm{z}) &\quad \text{if\ } x^{j} < z^{j},
		\end{cases}
\end{align} 
Under the same assumption over the highest coordinate, we find the second part of the generator $\Lan$ applied to $g_{\bm{z}}$ to be equal to
\begin{align}\label{2_g_z_nabla}
-\frac{1}{\alpha}\langle \bm{x},\nabla g_{\bm{z}}(\bm{x}) \rangle = 
		\begin{cases}
		-F_{\bm{Z}}(\bm{z})^{1 - \big(\frac{z^{j}}{x^{j}}\big)^{\alpha}} &\quad \text{if\ } x^{j} > z^{j}\\
		0 &\quad \text{if\ } x^{j} < z^{j}.
		\end{cases}
\end{align}
Summing the last two identities yields the desired result.

\end{proof}
\subsection{The case of the smooth Wasserstein distance}

We give the main properties of the Stein solution $g_{h}$ when $h \in \cilipone$.

\begin{proposition}\label{2_Lip}  
Let $\alpha>1$ and $h \in \cilipone$.

\begin{enumerate} 
\item  $\Pan_{t}h$ also belongs to $\cilipone$, and for all $j \in \lbra 1,d \rbra$:
\begin{align}\label{2_Pan_partl}
\vert \partial_{j}\Pan_{t}h(\bm{x}) \vert \leq e^{-\frac{t}{\alpha}}\prob(Z \leq \gamma_{t}^{-1/\alpha}x^{j}),\ t \geq 0,\ \bm{x} \in \R_{+}^{d}
\end{align}
with $Z\sim \mathcal{F}(\alpha)$. In particular the gradient of $\Pan_{t}h$ satisfies
\begin{align}\label{2_Pan_Lip}
\Vert\nabla \Pan_{t}h(\bm{x})\Vert_{1} \leq e^{-\frac{t}{\alpha}}\Vert \bm{p}_{t}(\bm{x})\Vert_{1}, 
\end{align}
where $\bm{p}_{t}(\bm{x}) \coloneq (p_{t}(x^{1}),\dots,p_{t}(x^{d})) = \big(\prob(Z \leq \gamma_{t}^{-1/\alpha}x^{1}),\dots,\prob(Z \leq \gamma_{t}^{-1/\alpha}x^{d})\big)$. As a result, $\Pan_{t}h$ is $de^{-\frac{t}{\alpha}}$-Lipschitz.

\item Assume that $h \in \mathrm{Lip}^{[2]}_{1}(\R_{+}^{d})$. Then $\partial_{j}\Pan_{t}h$ is $C_{1, \alpha}(t)$-Lipschitz for all $j \in \lbra 1,d \rbra$ and $t > 0$, with
\[
C_{1, \alpha}(t) \coloneq (d-1)e^{-t/\alpha} + \alpha\gamma_{t}^{-1/\alpha}.
\]

\end{enumerate}
\end{proposition}
 
\begin{proof}

1. We have seen already that $\Pan_{t}h$ is differentiable with respect to each argument, with $j$-th partial derivative equal to
\[
\partial_{j}\Pan_{t}h(\bm{x}) = e^{-\frac{t}{\alpha}}\esp\big[(\partial_{j}h)\big(e^{-\frac{t}{\alpha}}\bm{x} \oplus (1-e^{-t})^{\frac{1}{\alpha}}\bm{Z}\big)\ind_{\lb x^{j} \geq \gamma_{t}^{1/\alpha}Z^{j}\rb}\big]
\]
The fact that $\partial_{j}\Pan_{t}h(\bm{x})$ is continuous comes after noticing that $Z^{j}$ is diffuse, so that $\prob_{Z^{j}}$-a.s.
\[
\ind_{\lb x^{j} > \gamma_{t}^{1/\alpha}Z^{j}\rb} = \ind_{\lb x^{j} \geq \gamma_{t}^{1/\alpha}Z^{j}\rb} 
\]
Inequality \eqref{2_Pan_Lip} is a straightforward consequence of \eqref{2_Pan_partl_bnd}. 

2. The partial derivative $\partial_{j}\Pan_{t}h$ is $e^{-t/\alpha}$-Lipschitz with respect to $x^{k}$ for $k\neq j$ as a composition of the 1-Lipschitz function $\partial_{j}h$ and the $e^{-t/\alpha}$-Lipschitz function $x^{k} \mapsto e^{-t/\alpha}\bm{x}\oplus (1-e^{-t})^{1/\alpha}\bm{Z}$. To deal with the case $k=j$, we will assume $k = j = 1$. For all $\bm{x} \in \R_{+}^{d}$ and $\bm{y} = (y^{1}, x^{2},\dots, x^{d})$ we write
\begin{align*}
\partial_{1}\Pan_{t}h(\bm{x}) &- \partial_{1}\Pan_{t}h(\bm{y})\\
		&= \esp\big[\big(\partial_{1}h(\bm{Z}_{\bm{x}}) - \partial_{1}h(\bm{Z}_{\bm{y}})\big)\ind_{\lb x^{1} \geq \gamma_{t}^{1/\alpha}Z^{1}\rb}\big] + \esp\big[\partial_{1}h(Z_{\bm{y}})\big(\ind_{\lb x^{1} \geq \gamma_{t}^{1/\alpha}Z^{1}\rb} - \ind_{\lb y^{1} \geq \gamma_{t}^{1/\alpha}Z^{1}\rb}\big)\big].
\end{align*}
with $\bm{Z}_{\bm{x}} \coloneq e^{-t/\alpha}\bm{x} \oplus (1-e^{-t})^{1/\alpha}\bm{Z}$.
%&\leq e^{-\frac{t}{\alpha}}\vert x^{1} - y^{1}\vert + \esp\big[\vert\ind_{\lb x^{1} \geq \gamma_{t}^{1/\alpha}Z^{1}\rb} - \ind_{\lb y^{1} \geq \gamma_{t}^{1/\alpha}Z^{1}\rb}\vert\big]\\
%		&= e^{-\frac{t}{\alpha}}\vert x^{1} - y^{1}\vert + \big\vert F_{Z}\big(\gamma_{t}^{-1/\alpha}x^{1}\big) - F_{Z}\big(\gamma_{t}^{-1/\alpha}y^{1}\big)\big\vert.
The last identity comes from the monotony of $x \mapsto \ind_{\lb x \geq \gamma_{t}^{1/\alpha}Z^{1}\rb}$ and the fact that $F_{Z}$ is $\alpha$-Lipschitz.
\end{proof}

The next bound on the derivative of $g_{h}$ is helpful in applications. 

\begin{proposition}
Let $h \in \cilipone$. 
\begin{enumerate}
\item The associated Stein solution $g_{h}$ is $\alpha$-Lipschitz and its derivative satisfies:
\begin{align}\label{2_g_h_inqlty}
\vert \partial_{j}g_{h}(\bm{x}) \vert \leq \min\big(\alpha, (x^{j})^{\alpha}\big),\ \bm{x} \in \R_{+}^{d}.
\end{align}

\item Consequently, $\Dan g_{h}(\bm{x})$ is well-defined and finite for every $\bm{x} \in Eo$.

\item Furthermore, assume that $h \in \mathrm{Lip}^{[2]}_{1}(\R_{+}^{d})$. Then $\partial_{j}g_{h}$ is $C_{2, \alpha}$-Lipschitz, with 
\[
C_{2, \alpha} \coloneq \int_{0}^{\infty}C_{1, \alpha}(t)\dif t,
\]
and $C_{1, \alpha}(t)$ was defined in proposition \ref{2_Lip}.
\end{enumerate}
\end{proposition}

\begin{proof}

1. This is an immediate consequence of \eqref{2_Pan_partl} and of the definition of $g_{h}(\bm{x})$, the inversion between the partial derivative and the integral being justified as in theorem \ref{2_Ln-1}.

2. We have, again thanks to \eqref{2_Pan_partl}:
\begin{align*}
\int_{\Epol}\big\vert \big\langle r\bm{v}, \nabla g_{h}(\bm{x}) \big\rangle_{\bm{x}} \big\vert \dif \rho_{\alpha}(r) \dif\nu_{\alpha}(\bm{v})
		&= \sum_{j=1}^{d}\int_{\lb rv^{j} \geq x^{j}\rb} rv^{j} \vert\partial_{j}g_{h}(\bm{x})\vert \dif\rho_{\alpha}(r) \dif\nu_{\alpha}(\bm{v})\\
		&\leq \sum_{j=1}^{d}\int_{\Epol} rv^{j}\min(\alpha, r^{\alpha}) \frac{\alpha}{r^{\alpha + 1}}\dif r \dif\nu_{\alpha}(\bm{v})\\
		&= \sum_{j=1}^{d}\int_{\S_{+}^{d-1}}v^{j}\dif \nu(\bm{u})\int_{0}^{\infty}\min(\alpha, r^{\alpha}) \frac{\alpha}{r^{\alpha}}\dif r < +\infty.
\end{align*}
The integral over $\S^{d-1}_{+}$ equals $1$, thanks to the moment constraint relation \eqref{prelim_moment_constraints} satisfied by $\nu$, and the integral over $\R_{+}$ is finite because $\alpha > 1$. 

3. Recall the notations of proposition \ref{2_Lip}. We have seen that $\partial_{j}\Pan_{t}$ is $C_{1, \alpha}(t)$-differentiable. To justify the permutation between the partial derivative $\partial_{1}$ (say) and the integral sign, we write  
\begin{align*}
\Big\vert\partial_{1}\Pan_{t}h(\bm{x})\esp\Big[\partial_{1}h\big(\bm{Z}(\bm{x})\big)\ind_{\lb x^{1} \geq \gamma_{t}^{1/\alpha}Z^{1}\rb}\Big]\Big\vert &\leq \prob\big(Z \leq \gamma_{t}^{-1/\alpha}x^{1}\big) \leq \prob\big(Z \leq \gamma_{t}^{-1/\alpha}b\big) = e^{-\frac{\gamma_{t}}{b^{\alpha}}}, 
\end{align*}
where we assumed that $x^{1} \in [a, b]$, for $a \leq b$ two non-negative numbers and $Z \sim \mathcal{F}(\alpha)$. This last function is integrable with respect to $t$, so we can write:
\begin{align*}
\vert\partial_{1}g_{h}(\bm{x})\vert &\leq \int_{0}^{\infty}\vert\partial_{1}\Pan_{t}h^{*}(\bm{x})\vert\dif t\\
		&\leq \int_{0}^{\infty}C_{1, \alpha}(t)\dif t\\
		&= \int_{0}^{\infty}(d-1)e^{-\frac{t}{\alpha}} + \alpha\gamma_{t}^{-1/\alpha} \dif t < +\infty,
\end{align*}
because $\gamma^{-1/\alpha}_{t} \sim t^{-1/\alpha}$ as $t$ goes to $0^{+}$, and $1/\alpha \in (0,1)$.
\end{proof}

\subsection{The case of the Kolmogorov distance}

In this section we apply Stein's method to the target distribution $\mathcal{MS}(\alpha, \nu)$, for some $\alpha > 0$ and $\nu$ a finite measure on $\S_{+}^{d-1}$ satisfying the moment constraints \eqref{prelim_moment_constraints}.

\begin{proposition}\label{2_prop_g_z}
Let $\bm{z} \in \Eos$. We have the following:
\begin{enumerate}

\item $g_{\bm{z}}$ is constant over $[\bm{0}, \bm{z}]$:
\[
g_{\bm{z}}(\bm{x}) = -F_{\bm{Z}}(\bm{z})\int_{0}^{\infty}\big(F_{\bm{Z}}(\bm{z})^{-e^{-t}} - 1\big)\dif t,\ \bm{x} \in [0,\bm{z}].
\]
Furthermore we have the equivalent 
\[ 
g_{\bm{z}}(\bm{x}) \underset{x^{j} \to +\bm{\infty}}{\sim} \alpha F_{\bm{Z}}(\bm{z})\log x^{j},\ j \in \lbra 1,d \rbra.
\]
\item Assume there is only one $j_{0} \in \lbra 1,d \rbra$ such that
\[
\frac{x^{j_{0}}}{z^{j_{0}}} \geq \frac{x^{j}}{z^{j}},\ \forall j \in \lbra 1,d \rbra.
\]
Then $g_{\bm{z}}$ is continuously differentiable over $\R_{+}^{d}$ and one has:
\begin{align}\label{2_g_z_diff_2}
\partial_{j}g_{\bm{z}}(\bm{x}) = 
\begin{cases}
\frac{\alpha}{x^{j}}F_{\bm{Z}}(\bm{z})^{1 - \big(\frac{z^{j}}{x^{j}}\big)^{\alpha}}\quad & \text{if\ } j = j_{0}\\
0\quad & \text{otherwise}. 
\end{cases}
\end{align}
%Consequently, $\sup_{\bm{z} \in \Eo}\Da g_{\bm{z}}(\bm{0})$ is finite.

\item $\partial_{j}g_{\bm{z}}$ satisfies the inequality
\begin{align}\label{2_g_z_inqlty}
\vert \partial_{j}g_{\bm{z}}(\bm{x}) \vert \leq \frac{\alpha}{x^{j}}\ind_{\lb x^{j} > z^{j}\rb},\ j \in \lbra 1,d \rbra.
\end{align}
Thus it also satisfies for all $\bm{x},\bm{y} \in \R_{+}^{d}$:
\begin{align}\label{2_g_z_log}
\alpha^{-1}\vert g_{\bm{z}}(\bm{x}) - g_{\bm{z}}(\bm{y}) \vert \leq \Vert \log (\bm{x}\oplus \bm{z}) - \log (\bm{y}\oplus \bm{z}) \Vert_{1}\odot \Big(\frac{1}{z^{*}}\Vert \bm{x}\oplus \bm{z} - \bm{y}\oplus \bm{z}\Vert_{1}\Big),
\end{align}
where $z^{*}$ is the smallest non-null coordinate of $\bm{z}$. In particular, $g_{\bm{z}}$ is $\alpha/z^{*}$-Lipschitz.
\end{enumerate}
\end{proposition}

\begin{proof}

1. Identity \eqref{2_psdo-inv_Kol} makes obvious the fact that $g_{\bm{z}}$ is constant on $[\bm{0}, \bm{z}]$. As for the equivalent, the term inside the integral satisfies
\[
F_{\bm{Z}}(\bm{z})^{-e^{-t}} - 1 \underset{t \to +\infty}{\sim} F_{\bm{Z}}(\bm{z})^{-1}e^{-t}
\]
so that the integral part is equivalent to $F_{\bm{Z}}(\bm{z})^{-1}z^{j}/x^{j}$ when $x^{j}$ goes to infinity. 

2. The assumption on $j_{0}$ implies that $x^{j_{0}}$ is greater than $z^{j_{0}}$, so that the positive parts in \eqref{2_psdo-inv_Kol} vanish. Differentiating that identity yields the announced result.

3. The first inequality is a straightforward consequence of the previous point. The second one is obtained by integrating directly the inequality we just proved for every $j \in \lbra 1,d \rbra$, as well as by bounding $1/x^{j}$ by $1/z^{*}$ first.
\end{proof}

\section{Applications}

\subsection{Distance to a max-stable random vector}

In this subsection, the reference norm $\Vert\cdot\Vert$ will be any norm on $\R^{d}$ such that $\bm{x} = 1$ implies that every coordinate of $x^{j}$ is less than 1, \textit{e.g.} any $\ell^{p}$-norm. This choice is arbitrary and serves only to make certain arguments easier to state.
 
\subsubsection*{Distance between max-stable random vectors}

Let $\bm{Z}_{1}$ and $\bm{Z}_{2}$ be two max-stable random vectors with distribution $\mathcal{MS}(\alpha_{i}, \nu_{i})$ respectively, with $\alpha_{i} > 0$ and $i \in \lb 1,2 \rb$. We wish to get an idea about how close the distributions of $\bm{Z}_{1}$ and $\bm{Z}_{2}$ are, and express this closeness in terms of $\alpha_{i}$ and $\nu_{i}$. 

We start by assuming that $\nu_{1} = \nu_{2}$, \textit{i.e.} $\bm{Z}_{1}$ and $\bm{Z}_{2}$ have the same angular measure $\nu$ but different stability indices $\alpha_{1}, \alpha_{2}$. 

\begin{proposition}\label{3_dist_alpha}
Let $\alpha_{1}, \alpha_{2}$ be two positive numbers with $\alpha_{1} < \alpha_{2}$, and $\nu$ a finite measure on $\S_{+}^{d-1}$ satisfying the moment constraints \eqref{prelim_moment_constraints}. Let $Z_{i} \sim \mathcal{MS}(\alpha_{i}, \nu)$ for $i = 1,2$.
\begin{itemize}
\item[-] There exists a constant $C^{\mathrm{K}}_{\alpha_{1}, \alpha_{2}, \nu}$ such that:
\begin{align*}
\mathrm{d}_{\mathrm{K}}(\bm{Z}_{1}, \bm{Z}_{2}) \leq C^{\mathrm{K}}_{\alpha_{1}, \alpha_{2}, \nu}\Big(\frac{1}{\alpha_{1}} + \frac{1}{\alpha_{2}}\Big)\vert \alpha_{1} - \alpha_{2}\vert.
\end{align*} 
Let $Z_{1} \sim \mathcal{F}(\alpha_{1})$ and $\mu_{\alpha_{1}}$ the exponent measure of $\bm{Z}_{1}$. A possible choice for the constant $C^{\mathrm{K}}_{\alpha_{1}, \alpha_{2}, \nu_{1}}$ is:
\begin{align}\label{3_C_an_K}
C^{\mathrm{K}}_{\alpha_{1}, \alpha_{2}, \nu} = d + \sum_{j=1}^{d}\int_{\Eo}\vert\log y^{j}\vert\prob\big(Z_{1}^{\alpha_{2}/\alpha_{1}}\odot Z_{1} \leq y^{j}\big)\dif\mu_{\alpha_{1}}(\bm{y}).
\end{align}

\item[-] If furthermore both $\alpha_{1}$ and $\alpha_{2}$ are greater than $1$, we have a constant $C^{\mathrm{W}}_{\alpha_{1}, \alpha_{2}, \nu_{1}} > 0$ such that:
\begin{align*}
\mathrm{d}_{\mathrm{W}}(\bm{Z}_{1}, \bm{Z}_{2}) \leq C^{\mathrm{W}}_{\alpha_{1}, \alpha_{2}, \nu}\Big(\frac{1}{\alpha_{1}} + \frac{1}{\alpha_{2}}\Big)\vert \alpha_{1} - \alpha_{2}\vert.
\end{align*} 
A possible choice for the constant $C_{\alpha_{1}, \alpha_{2}, \nu}$ is:
\begin{align}\label{3_C_an_W}
C^{\mathrm{W}}_{\alpha_{1}, \alpha_{2}, \nu} = d\Gamma\Big(1 - \frac{1}{\alpha_{1}}\Big) + \sum_{j=1}^{d}\int_{\Eo}(1 \oplus y^{j})\vert \log y^{j} \vert\prob\big(Z_{1}^{\alpha_{2}/\alpha_{1}}\odot Z_{1} \leq y^{j}\big)\dif\mu_{\alpha_{1}}(\bm{y}).
\end{align}
\end{itemize}
\end{proposition}

\begin{proof}

- Set $h^{*} = h - \esp[h(\bm{Z}_{2})]$ for any $\prob_{\bm{Z}_{2}}$-integrable $h$, and take $g_{\bm{z}} = \mathscr{L}_{\alpha_{2}, \nu}^{-1}h^{*}_{\bm{z}}$. We have:
\begin{align*}
\vert F_{\bm{Z}_{1}}(\bm{z}) - F_{\bm{Z}_{2}}(\bm{z}) \vert &= \vert \esp[h_{\bm{z}}(\bm{Z}_{1})] - \esp[h_{\bm{z}}(\bm{Z}_{2})]\vert\\
		&= \vert \esp[h^{*}_{\bm{z}}(\bm{Z}_{1})]\vert\\
		&= \vert \esp[\mathscr{L}_{\alpha_{2}, \nu}g_{\bm{z}}(\bm{Z}_{1})]\vert.
\end{align*}
Thanks to inequality \eqref{2_g_z_inqlty}, we see that $\langle \bm{Z}_{1}, \nabla g_{\bm{z}}(\bm{Z}_{1}) \rangle$ and $\mathbf{D}_{\alpha_{2}, \nu}g_{\bm{z}}(\bm{Z}_{1})$ have finite expectations. Furthermore, the fact that $\esp[\mathscr{L}_{\alpha_{1},\nu}g_{\bm{z}}(\bm{Z}_{1})] = 0$ implies that:
\[
\esp\big[\mathbf{D}_{\alpha_{1}, \nu}g_{\bm{z}}(\bm{Z}_{1})\big] = \frac{1}{\alpha_{1}}\esp\big[\langle \bm{Z}_{1}, \nabla g_{\bm{z}}(\bm{Z}_{1})\rangle\big],
\]
so that
\begin{align*}\label{3_dist_a1_a2}
\alpha_{2}\esp\big[\mathscr{L}_{\alpha_{2}, \nu}g_{\bm{z}}(\bm{Z}_{1})\big] &= \alpha_{2}\esp\big[\mathbf{D}_{\alpha_{2}, \nu}g_{\bm{z}}(\bm{Z}_{1})\big] - \esp\big[\langle \bm{Z}_{1}, \nabla g_{\bm{z}}(\bm{Z}_{1})\rangle\big]\\
		&= \alpha_{2}\esp\big[\mathbf{D}_{\alpha_{2}, \nu}g_{\bm{z}}(\bm{Z}_{1}) - \mathbf{D}_{\alpha_{1}, \nu}g_{\bm{z}}(\bm{Z}_{1})\big] + \alpha_{2}\esp[\mathbf{D}_{\alpha_{1}, \nu}g_{\bm{z}}(\bm{Z}_{1})] - \esp\big[\langle \bm{Z}_{1}, \nabla g_{\bm{z}}(\bm{Z}_{1})\rangle\big]\\
		&= \alpha_{2}\esp\big[\mathbf{D}_{\alpha_{2}, \nu}g_{\bm{z}}(\bm{Z}_{1}) - \mathbf{D}_{\alpha_{1}, \nu}g_{\bm{z}}(\bm{Z}_{1})\big] + \Big(\frac{\alpha_{2}}{\alpha_{1}} - 1\Big)\esp\big[\langle \bm{Z}_{1}, \nabla g_{\bm{z}}(\bm{Z}_{1})\rangle\big]. \numberthis
\end{align*}
The second term is bounded by $d\alpha_{2}\vert 1 - \alpha_{2}\alpha_{1}^{-1} \vert$, for
\[
\esp\big[\langle \bm{Z}_{1}, \nabla g_{\bm{z}}(\bm{Z}_{1})\rangle\big] = \sum_{j=1}^{d}\esp\big[Z^{j}_{1}\partial_{j}g_{\bm{z}}(\bm{Z}_{1})\big] \leq d\alpha_{2}, 
\]
thanks to inequality \eqref{2_g_z_inqlty}. As for the first one, comparing directly the operators $\mathbf{D}_{\alpha_{2}, \nu}$ and $\mathbf{D}_{\alpha_{1}, \nu}$ seems difficult. Instead, we start by giving a more convenient expression of $\mathbf{D}_{\alpha_{2}, \nu}g_{\bm{z}}(\bm{Z}_{1})$. As in \cite{Costaceque24_mlti}, we denote by $T_{\alpha}$ for any $\alpha > 0$ the application:
\[
T_{\alpha}f(\bm{x}) = f(\bm{x}^{\alpha}).
\]
Since $\Pan = T_{\alpha}\Pn T_{\alpha^{-1}}$, we can write:
\begin{align*}
\esp[\mathbf{D}_{\alpha_{2}, \nu}g_{\bm{z}}(\bm{Z}_{1})] &= \esp\Big[T_{\frac{\alpha_{2}}{\alpha_{1}}}\mathbf{D}_{\alpha_{1}, \nu}T_{\frac{\alpha_{1}}{\alpha_{2}}}g_{\bm{z}}(\bm{Z}_{1})\Big]\\
		&= \esp\Big[\mathbf{D}_{\alpha_{1}, \nu}T_{\frac{\alpha_{1}}{\alpha_{2}}}g_{\bm{z}}\big(\bm{Z}_{1}^{\frac{\alpha_{2}}{\alpha_{1}}}\big)\Big]\\
		&= \esp\big[\mathbf{D}_{\alpha_{1}, \nu}T_{\beta^{-1}}g_{\bm{z}}(\bm{Z}_{1}^{\beta})\big],
\end{align*}
where we set $\beta \coloneq \alpha_{2}/\alpha_{1}$. Let $\mu_{\alpha_{1}}$ be the exponent measure of $\bm{Z}_{1}$. From what precedes we get
\begin{align*}
\esp\big[\mathbf{D}_{\alpha_{2}, \nu}g_{\bm{z}}(\bm{Z}_{1}) - \mathbf{D}_{\alpha_{1}, \nu}g_{\bm{z}}(\bm{Z}_{1})\big]	&= \esp\big[\mathbf{D}_{\alpha_{1}, \nu}T_{\beta^{-1}}g_{\bm{z}}(\bm{Z}_{1}^{\beta}) - \mathbf{D}_{\alpha_{1}, \nu}g_{\bm{z}}(\bm{Z}_{1})\big]\\
		&= \int_{\Eo}\esp\big[g_{\bm{z}}\big(\bm{Z}_{1} \oplus \bm{y}^{1/\beta}\big) - g_{\bm{z}}(\bm{Z}_{1}) - g_{\bm{z}}\big(\bm{Z}_{1} \oplus \bm{y}\big) + g_{\bm{z}}(\bm{Z}_{1})\big]\dif\mu_{\alpha_{1}}(\bm{y})\\
		&= \int_{\Eo}\esp\big[g_{\bm{z}}\big(\bm{Z}_{1} \oplus \bm{y}^{1/\beta}\big) - g_{\bm{z}}\big(\bm{Z}_{1} \oplus \bm{y}\big)\big]\dif\mu_{\alpha_{1}}(\bm{y}).
\end{align*}
The triangle inequality and bound \ref{2_g_z_log} yield:
\begin{align*}
\big\vert\esp\big[\mathbf{D}_{\alpha_{2}, \nu}g_{\bm{z}}(\bm{Z}_{1}) - \mathbf{D}_{\alpha_{1}, \nu}g_{\bm{z}}(\bm{Z}_{1})\big]\big\vert &\leq \int_{\Eo}\esp\big[\big\vert g_{\bm{z}}\big(\bm{Z}_{1} \oplus \bm{y}^{1/\beta}\big) - g_{\bm{z}}\big(\bm{Z}_{1} \oplus \bm{y}\big)\big\vert\big]\dif\mu_{\alpha_{1}}(\bm{y})\\
		&\leq \alpha_{2}\int_{\Eo}\esp\big[\big\Vert \log\big(\bm{Z}_{1}\oplus \bm{y}^{1/\beta}\big) - \log\big(\bm{Z}_{1}\oplus \bm{y}\big)\big\Vert_{1}\big]\dif\mu_{\alpha_{1}}(\bm{y})	\\
		&= \alpha_{2}\sum_{j=1}^{d}\int_{\Eo}\esp\big[\big\vert \log\big(Z_{1}\oplus (y^{j})^{1/\beta}\big) - \log\big(Z_{1}\oplus y^{j}\big)\big\vert\big]\dif\mu_{\alpha_{1}}(\bm{y})\\
		&= \alpha_{2}\sum_{j=1}^{d}\int_{\Eo}\esp\big[\big\vert \log\big(Z_{1}\oplus (y^{j})^{1/\beta}\big) - \log\big(Z_{1}\oplus y^{j}\big)\big\vert\ind_{\lb y^{j} \geq Z_{1}^{\beta}\odot Z_{1}\rb}\big]\dif\mu_{\alpha_{1}}(\bm{y}),		
\end{align*}
where $Z_{1} \sim \mathcal{F}(\alpha_{1})$. The indicator function in the last equality comes from noticing that if $y^{j} \leq Z_{1}$ and $y^{j} \leq Z_{1}^{\beta}$ at the same time, then the integrand vanishes. Finally, because the logarithm is a non-decreasing function, we get that $\log(x\oplus y) = \log x \oplus \log y$. Furthermore, $x\mapsto c\oplus x$ is $1$-Lipschitz. Using this result with $c = \log Z_{1}$, we find:
\begin{align*}
\big\vert\esp\big[\mathbf{D}_{\alpha_{2}, \nu}g_{\bm{z}}(\bm{Z}_{1}) - \mathbf{D}_{\alpha_{1}, \nu}g_{\bm{z}}(\bm{Z}_{1})\big]\big\vert &\leq \alpha_{2}\sum_{j=1}^{d}\int_{\Eo}\vert \beta^{-1}\log y^{j} - \log y^{j}\vert\prob\big(Z_{1}^{\beta}\odot Z_{1} \leq y^{j}\big)\dif\mu_{\alpha_{1}}(\bm{y})\\
		&= \vert \alpha_{2} - \alpha_{1}\vert \sum_{j=1}^{d}\int_{\Eo}\vert\log y^{j}\vert\prob\big(Z_{1}^{\beta}\odot Z_{1} \leq y^{j}\big)\dif\mu_{\alpha_{1}}(\bm{y}).
\end{align*}
This last integral is bounded since the probability inside vanishes at exponential speed when $y^{j}$ goes to $0$. By dividing this bound by $\alpha_{2}$, we obtain the desired result.

%This term in the integral also writes for each $j \in \lbra 1,d \rbra$:
%\begin{align*}
%\int_{\Eo}\vert\log y^{j}\vert\prob\big(Z_{1}^{\beta}\odot Z_{1} \leq y^{j}\big)\dif\mu_{\alpha_{1}}(\bm{y}) &= \int_{\lb y^{j} \leq 1 \rb}\vert\log y^{j}\vert\prob\big(Z_{1}^{\beta}\odot Z_{1} \leq y^{j}\big)\dif\mu_{\alpha_{1}}(\bm{y})\\
%		&\quad + \int_{\lb y^{j} \geq 1 \rb}\vert\log y^{j}\vert\prob\big(Z_{1}^{\beta}\odot Z_{1} \leq y^{j}\big)\dif\mu_{\alpha_{1}}(\bm{y})\\
%		&\leq \int_{\lb y^{j} \leq 1 \rb}\vert\log y^{j}\vert\prob\big(Z_{1}^{\beta} \leq y^{j}\big)\dif\mu_{\alpha_{1}}(\bm{y})\\
%		&\quad + \int_{\lb y^{j} \geq 1 \rb}\log y^{j}\dif\mu_{\alpha_{1}}(\bm{y}).
%\end{align*}

- We start from equality \eqref{3_dist_a1_a2}, but with $g_{h} = \mathscr{L}_{\alpha_{2},\nu}^{-1}h^{*}$ instead of $g_{\bm{z}}$, where $h$ is in $\cilipone$. As in the previous point, we have two terms to bound:
\begin{align*}
\big\vert\esp\big[\langle \bm{Z}_{1}, \nabla g_{\bm{h}}(\bm{Z}_{1})\rangle\big]\big\vert &\leq \sum_{j=1}^{d}\esp\big[Z_{1}^{j}\vert\partial_{j}g_{\bm{h}}(\bm{Z}_{1})\vert\big]\\
	&\leq \alpha_{2}\sum_{j=1}^{d}\esp[Z_{1}^{j}]\\
	&= d\alpha_{2}\Gamma\Big(1 - \frac{1}{\alpha_{1}}\Big),
\end{align*}
this inequality being a consequence of \eqref{2_g_h_inqlty}. Concerning the second term, we use the Lipschitz property of $g_{h}$.
%relying on the fact that both $\alpha_{1}$ and $\alpha_{2}$ are greater than $1$. 
\begin{align*}
\big\vert\esp\big[\mathbf{D}_{\alpha_{2}, \nu}g_{h}(\bm{Z}_{1}) - \mathbf{D}_{\alpha_{1}, \nu}g_{h}(\bm{Z}_{1})\big]\big\vert &\leq \int_{\Eo}\esp\Big[\big\vert g_{h}\big(\bm{Z}_{1} \oplus \bm{y}^{1/\beta}\big) - g_{h}\big(\bm{Z}_{1} \oplus \bm{y}\big)\big\vert\Big]\dif\mu_{\alpha_{1}}(\bm{y})\\
		&\leq \alpha_{2}\int_{\Eo}\esp\Big[\big\Vert \bm{Z}_{1} \oplus \bm{y}^{1/\beta} - \bm{Z}_{1} \oplus \bm{y}\big\Vert_{1}\Big]\dif\mu_{\alpha_{1}}(\bm{y})\\
		&= \alpha_{2}\sum_{j=1}^{d}\int_{\Eo}\esp\Big[\big\vert Z^{j}_{1} \oplus (y^{j})^{1/\beta} - Z^{j}_{1} \oplus y^{j}\big\vert\ind_{\lb y^{j} \geq Z_{1}^{\beta}\odot Z_{1}\rb}\Big]\dif\mu_{\alpha_{1}}(\bm{y})\\
		&\leq \alpha_{2}\sum_{j=1}^{d}\int_{\Eo}\big\vert(y^{j})^{\beta^{-1}} - y^{j}\big\vert\prob\big(Z_{1}^{\beta}\odot Z_{1} \leq y^{j}\big)\dif\mu_{\alpha_{1}}(\bm{y})\\
		&= \alpha_{2}\sum_{j=1}^{d}\int_{\Eo}y^{j}\big\vert(y^{j})^{\beta^{-1}-1} - 1\big\vert\prob\big(Z_{1}^{\beta}\odot Z_{1} \leq y^{j}\big)\dif\mu_{\alpha_{1}}(\bm{y}).
\end{align*}
This quantity is finite because $\alpha_{1}>1$ and $\mu_{\alpha_{1}}$ is the image measure of $\rho_{1}\otimes \nu$ by $(r, \bm{u}) \mapsto (r\bm{u})^{1/\alpha_{1}}$, with $\nu$ a finite measure over the compact $\S_{+}^{d-1}$. Recall that $\alpha_{1} < \alpha_{2}$, so that $\beta^{-1} < 1$, and thus $1-\beta^{-1} \in (0,1)$. We need to study the function $f : \delta \mapsto y^{-\delta} = \exp(-\delta\log y)$ over $(0,1)$ for all $y \in \R_{+}^{*}$. If $y \in [1,+\infty)$, we have
\[
\vert y^{-\delta} - 1 \vert = 1 - \exp(-\delta\log y) \leq \delta\log y.
\]
And if $y \in (0,1]$, then $\vert f'(\delta) \vert = y^{-\delta}\vert \log y \vert \leq y^{-1}\vert\log y \vert$. Thus $f$ is $(1\oplus y^{-1})\vert \log y \vert$-Lipschitz over $\R_{+}^{*}$. Finally:
\begin{align*}
\big\vert\esp\big[\mathbf{D}_{\alpha_{2}, \nu}g_{h}(\bm{Z}_{1}) - \mathbf{D}_{\alpha_{1}, \nu}g_{h}(\bm{Z}_{1})\big]\big\vert &= \alpha_{2}\sum_{j=1}^{d}\int_{\Eo}y^{j}\big\vert(y^{j})^{\beta^{-1}-1} - 1\big\vert\prob\big(Z_{1}^{\beta}\odot Z_{1} \leq y^{j}\big)\dif\mu_{\alpha_{1}}(\bm{y})\\
		&\leq \vert\alpha_{1} - \alpha_{2}\vert\sum_{j=1}^{d}\int_{\Eo}(1 \oplus y^{j})\vert \log y^{j} \vert\prob\big(Z_{1}^{\beta}\odot Z_{1} \leq y^{j}\big)\dif\mu_{\alpha_{1}}(\bm{y}),
\end{align*}
and this last integral is finite.
\end{proof}

Next we suppose the stability index is now the same $\alpha$ for both $\bm{Z}_{1}$ and $\bm{Z}_{2}$ but that they differ by their angular measures, $\nu_{1}$ and $\nu_{2}$ respectively.

\begin{proposition}\label{3_dist_nu}
Let $\alpha$ be a positive number and $\nu_{1}, \nu_{2}$ two finite measures on $\S_{+}^{d-1}$ satisfying the moment constraints \eqref{prelim_moment_constraints}.  Let also $\bm{Z}_{i} \sim \mathcal{MS}(\alpha, \nu_{i})$ for $i = 1,2$.
\begin{itemize}

\item[-] The following inequality holds:
\begin{align*}
\mathrm{d}_{\mathrm{K}}(\bm{Z}_{1}, \bm{Z}_{2}) \leq d\mathrm{d}_{\mathrm{TV}}(\nu_{1}, \nu_{2}).
\end{align*} 

\item[-] If furthermore $\alpha$ is greater than $1$, then 
\begin{align*}
\mathrm{d}_{\mathrm{W}}(\bm{Z}_{1}, \bm{Z}_{2}) \leq d\Gamma\Big(1 - \frac{1}{\alpha}\Big)\mathrm{d}_{\mathrm{TV}}(\nu_{1}, \nu_{2}).
\end{align*}
\end{itemize}
\end{proposition}

\begin{proof}

- As before, set $h^{*} = h - \esp[h(\bm{Z}_{2})]$ for any $\prob_{\bm{Z}_{2}}$-integrable $h$, and $g_{\bm{z}} = \mathscr{L}_{\alpha_{2}, \nu}^{-1}h^{*}_{\bm{z}}$. We have
\[
\vert F_{\bm{Z}_{1}}(\bm{z}) - F_{\bm{Z}_{2}}(\bm{z}) \vert = \vert \esp[\mathscr{L}_{\alpha, \nu_{2}}g_{\bm{z}}(\bm{Z}_{1})]\vert.
\]
We use once again that the term $\esp\big[\langle \bm{Z}_{1}, \nabla g_{\bm{z}}(\bm{Z}_{1})\rangle\big]$ is common to both $\mathscr{L}_{\alpha, \nu_{1}}$ and $\mathscr{L}_{\alpha, \nu_{2}}$:
\begin{align*}\label{3_dist_nu1_nu2}
\esp\big[\mathscr{L}_{\alpha_{2}, \nu}g_{\bm{z}}(\bm{Z}_{1})\big] &= \esp\big[\mathbf{D}_{\alpha, \nu_{2}}g_{\bm{z}}(\bm{Z}_{1})\big] - \frac{1}{\alpha}\esp\big[\langle \bm{Z}_{1}, \nabla g_{\bm{z}}(\bm{Z}_{1})\rangle\big]\\
		&= \esp\big[\mathbf{D}_{\alpha, \nu_{2}}g_{\bm{z}}(\bm{Z}_{1}) - \mathbf{D}_{\alpha, \nu_{1}}g_{\bm{z}}(\bm{Z}_{1})\big]\numberthis.
\end{align*}
This time, comparing $\mathbf{D}_{\alpha, \nu_{1}}$ and $\mathbf{D}_{\alpha, \nu_{2}}$ is easier:
\begin{align*}
\mathbf{D}_{\alpha, \nu_{2}}g_{\bm{z}}(\bm{x}) - \mathbf{D}_{\alpha, \nu_{1}}&g_{\bm{z}}(\bm{x})\\
		&= \frac{1}{\alpha}\int_{\Epol}\big\langle r\bm{u}^{1/\alpha},\nabla g_{\bm{z}}(\bm{x}\oplus r\bm{u}^{1/\alpha})\big\rangle_{\bm{x}}\dif\rho_{\alpha}(r) \big(\mathrm{d}\nu_{2}(\bm{u}) - \dif\nu_{1}(\bm{u})\big)\\
		&= \frac{1}{\alpha}\int_{\R_{+}^{*}}\int_{\S_{+}^{d-1}}\big\langle r\bm{u}^{1/\alpha},\nabla g_{\bm{z}}(\bm{x}\oplus r\bm{u}^{1/\alpha})\big\rangle_{\bm{x}}\big(f_{\nu_{2}}(\bm{u}) - f_{\nu_{1}}(\bm{u})\big)\dif \nu_{1,2}(\bm{u})\dif\rho_{\alpha}(r).
\end{align*}
for every $\bm{x} \in \R_{+}^{d}$. Here $\nu_{1,2} \coloneq \nu_{1} + \nu_{2}$, and $f_{\nu_{i}} = \dif\nu_{i}/\dif\nu_{1,2}$ for $i = 1,2$, the density function of $\nu_{i}$ with respect to to $\nu_{1,2}$. Now, let $Z$ be a random variable with Fréchet distribution $\mathcal{F}(\alpha)$. Replacing $\bm{x}$ by $\bm{Z}_{1}$ and taking expectations, we bound the scalar product in the integral by
\begin{align*}
\esp\Big[\big\vert \big\langle r\bm{u}^{1/\alpha},\nabla g_{\bm{z}}(\bm{Z}_{1}\oplus r\bm{u}^{1/\alpha})\big\rangle_{\bm{Z}_{1}} \big\vert\Big] &\leq \alpha\sum_{j=1}^{d}\prob\big(Z \leq r(u^{j})^{1/\alpha}\big) \leq d\alpha\prob(Z \leq r),
\end{align*}
because 
\[
r(u^{j})^{1/\alpha}\vert\partial_{j}g_{\bm{z}}(\bm{Z}_{1}\oplus r\bm{u}^{1/\alpha})\vert \leq \alpha,
\]
thanks to inequality \eqref{2_g_z_inqlty}. We have also used that $u^{j}$ is always in $[0,1]$, as $\bm{u} \in \S_{+}^{d-1}$, due to our initial assumption on the reference norm. Finally we find:
\begin{align*}
\vert F_{\bm{Z}_{1}}(\bm{z}) - F_{\bm{Z}_{2}}(\bm{z}) \vert &= \big\vert \esp\big[\mathbf{D}_{\alpha, \nu_{2}}g_{\bm{z}}(\bm{Z}_{1}) - \mathbf{D}_{\alpha, \nu_{1}}g_{\bm{z}}(\bm{Z}_{1})\big] \big\vert\\
		&\leq d\int_{\R_{+}^{*}}\int_{\S_{+}^{d-1}}\prob(Z \leq r)\big\vert f_{\nu_{2}}(\bm{u}) - f_{\nu_{1}}(\bm{u})\big\vert\dif \nu_{1,2}(\bm{u})\dif\rho_{\alpha}(r)\\
		&= d\int_{\R_{+}^{*}}e^{-\frac{1}{r^{\alpha}}}\frac{\alpha}{r^{\alpha + 1}}\dif r\int_{\S_{+}^{d-1}}\big\vert f_{\nu_{2}}(\bm{u}) - f_{\nu_{1}}(\bm{u})\big\vert\dif \nu_{1,2}(\bm{u})\\
		&= d\mathrm{d}_{\mathrm{TV}}(\nu_{1}, \nu_{2}).
\end{align*}

- We start from \eqref{3_dist_nu1_nu2}, with $g_{h} = \mathscr{L}_{\alpha, \nu_{2}}h^{*}$ instead of $g_{\bm{z}}$, where $h : \R_{+}^{d} \to \R$ belongs to $\cilipone$. Then, we know from inequality \eqref{2_g_h_inqlty} that:
\begin{align*}
\big\vert \big\langle r\bm{u}^{1/\alpha},\nabla g_{h}(\bm{Z}_{1}\oplus r\bm{u}^{1/\alpha})\big\rangle_{\bm{Z}_{1}} \big\vert &\leq \sum_{j=1}^{d}\min(\alpha , r^{\alpha}u^{j})r(u^{j})^{1/\alpha}\prob\big(Z \leq r(u^{j})^{1/\alpha}\big)\\
		&\leq d\alpha r\prob(Z \leq r).
\end{align*} 
Consequently we find:
\begin{align*}
\vert h(\bm{Z}_{1}) - h(\bm{Z}_{2}) \vert &\leq d\int_{\R_{+}^{*}}e^{-\frac{1}{r^{\alpha}}}\frac{\alpha}{r^{\alpha}}\dif r\int_{\S_{+}^{d-1}}\big\vert f_{\nu_{2}}(\bm{u}) - f_{\nu_{1}}(\bm{u})\big\vert\dif \nu_{1,2}(\bm{u})\\
		&= d\Gamma\Big(1 - \frac{1}{\alpha}\Big)\mathrm{d}_{\mathrm{TV}}(\nu_{1}, \nu_{2}).
\end{align*}
\end{proof}

We combine those two results in the following corollary.

\begin{corollary}
Let $\alpha_{1}, \alpha_{2}$ be two positive numbers with $\alpha_{1} < \alpha_{2}$, and $\nu_{1}, \nu_{2}$ two finite measures on $\S_{+}^{d-1}$ satisfying the moment constraints \eqref{prelim_moment_constraints}. Let $\bm{Z}_{i} \sim \mathcal{MS}(\alpha_{i}, \nu_{i})$, for $i = 1,2$.
\begin{itemize}

\item[-] By taking $C^{\mathrm{K}}_{\alpha_{1}, \alpha_{2}, \nu_{1}}$ as in \eqref{3_C_an_K}, the following inequality holds:
\begin{align*}
\mathrm{d}_{\mathrm{K}}(\bm{Z}_{1}, \bm{Z}_{2}) \leq C^{\mathrm{K}}_{\alpha_{1}, \alpha_{2}, \nu_{1}}\Big(\frac{1}{\alpha_{1}} + \frac{1}{\alpha_{2}}\Big)\vert \alpha_{1} - \alpha_{2}\vert + d\mathrm{d}_{\mathrm{TV}}(\nu_{1}, \nu_{2}).
\end{align*} 

\item[-] If furthermore $\alpha_{1}, \alpha_{2}$ are greater than $1$, then for $C^{\mathrm{W}}_{\alpha_{1}, \alpha_{2}, \nu_{1}}$ given by \eqref{3_C_an_W}, we have: 
\begin{align*}
\mathrm{d}_{\mathrm{W}}(\bm{Z}_{1}, \bm{Z}_{2}) \leq C^{\mathrm{W}}_{\alpha_{1}, \alpha_{2}, \nu_{1}}\Big(\frac{1}{\alpha_{1}} + \frac{1}{\alpha_{2}}\Big)\vert \alpha_{1} - \alpha_{2}\vert + d\Gamma\Big(1 - \frac{1}{\alpha_{2}}\Big)\mathrm{d}_{\mathrm{TV}}(\nu_{1}, \nu_{2}).
\end{align*}
\end{itemize}
\end{corollary}

\begin{proof}
Simply bound $\mathrm{d}_{\mathrm{K}}(\bm{Z}_{1}, \bm{Z}_{2})$ by $\mathrm{d}_{\mathrm{K}}(\bm{Z}_{1}, \bm{Z}_{3}) + \mathrm{d}_{\mathrm{K}}(\bm{Z}_{3}, \bm{Z}_{2})$, with $\bm{Z}_{3}\sim \mathcal{MS}(\alpha_{2}, \nu_{1})$ and use propositions \ref{3_dist_alpha} and \ref{3_dist_nu}.
\end{proof}

\subsection{Rate of convergence of the de Haan-LePage series}

Let $n\in \N^{*}$ and $\phi = ((r_{i},\bm{v}_{i}))_{i \geq 1}$ be a configuration in $\Eo$. We arrange $\phi$ in decreasing order with respect to the first coordinate, \textit{i.e.} $r_{1}\geq r_{2} \geq \dots $, and define
\[ 
\m_{n}(\phi) \coloneq \bigoplus_{i=1}^{n}r_{i}\bm{v}_{i}. 
\]
Take $\bm{Z}$ a max-stable random vector with distribution $\mathcal{MS}(\alpha, \nu)$. Define a measure $\nu_{\alpha}$ on $\S_{+}^{d-1}$ by setting
\[
\nu_{\alpha}(B) \coloneq \nu(B^{\alpha}),\ B \in \mathscr{B}(\S_{+}^{d-1}),
\]
where $B^{\alpha} = \lb x^{\alpha},\ x \in B\rb$. In other words, $\nu_{\alpha}$ is the image-measure of $\nu$ by $\bm{u} \mapsto \bm{u}^{1/\alpha}$. Let $\eta = ((r_{i}, \bm{v}_{i}))_{i \geq 1}$ be a marked Poisson process on $\Epol$ with intensity measure $\alpha r^{-(\alpha + 1)}\dif r\otimes \dif\nu_{\alpha}(\bm{v})$. We know from \cite{Davydov08} that when $n$ goes to infinity, $\m_{n}(\eta)$ converges to $\m(\eta)$ in distribution. Our goal is to determine an estimation of the speed at which this convergence occurs.

\begin{proposition}\label{3_rate_LePage_Wass}
Let $\eta = ((r_{i}, \bm{v}_{i}))_{i\geq 1}$ be a Poisson process on $\Epol$ with measure $\alpha r^{-(\alpha+1)}\dif r\otimes \dif \nu_{\alpha}(\bm{v})$, with $\alpha > 1$. Set 
\[
\bm{Z}_{n} \coloneq \m_{n}(\eta) = \bigoplus_{i=1}^{n}r_{i}\bm{v}_{i}.
\]
Then there exists a constant $C_{\alpha}>0$ depending only on $d$ and $\alpha$ such that
\[ 
\mathrm{d}_{\mathrm{W}}(\bm{Z}_{n}, \bm{Z}) \leq C_{\alpha}\frac{1}{n},\ n\geq 2. 
\]
\end{proposition}

\begin{proof}

We first check that $\langle \bm{Z}_{n}, \nabla g_{h}(\bm{Z}_{n})\rangle$ is integrable. Thanks to property \ref{2_g_h_inqlty}, we know that $g_{h}$ is $\alpha$-Lipschitz for $h$ is $1$-Lipschitz; therefore 
\[
\vert \langle \bm{Z}_{n}, \nabla g_{h}(\bm{Z}_{n})\rangle \vert \leq \alpha\sum_{j=1}^{d}Z_{n}^{j} \leq \alpha\sum_{j=1}^{d}Z^{j}
\]
because by definition $Z_{n}^{j}$ is dominated by $\bm{Z} = \m(\eta)$, whose coordinates are integrable since $\alpha > 1$. The integrability of $\Dan f(\bm{Z}_{n})$ has already been checked in proposition \ref{2_g_h_inqlty}.

To compare $\esp[\Dan g_{h}(\bm{Z}_{n})]$ and $\alpha\esp[\langle \bm{Z}_{n}, g_{h}(\bm{Z}_{n})\rangle]$, we apply the Campbell-Mecke formula to:
\[ 
\eta \mapsto \big\langle r\bm{v}, \nabla g_{h}\big(\m_{n}(\eta)\oplus r\bm{v}\big)\big\rangle_{\m_{n}(\eta)}, 
\]
so that we get 
\begin{align*}\label{3_LePage_3parts}
\esp\big[\Dan g_{h}\big(\m_{n}(\eta)\big)\big] &= \int_{\Epol} \esp\Big[\big\langle r\bm{v}, \nabla g_{h}\big(\m_{n}(\eta)\oplus r\bm{v}\big)\big\rangle_{\m_{n}(\eta)}\Big] \frac{\alpha}{r^{\alpha+1}}\dif r\dif\nu_{\alpha}(\bm{v})\\
		&= \esp\Big[\sum_{(r,\bm{v}) \in \eta} \big\langle r\bm{v}, \nabla g_{h}\big(\m_{n}(\eta - \delta_{(r,\bm{v})})\oplus r\bm{v}\big)\big\rangle_{\m_{n}(\eta - \delta_{(r,\bm{v})})}\Big]\\
		&= \esp\Big[\sum_{(r,\bm{v}) \in \eta} \big\langle r\bm{v}, \nabla g_{h}\big(\m_{n}(\eta_{r,\bm{v}}))\oplus r\bm{v}\big)\big\rangle_{\m_{n}(\eta_{r,\bm{v}})}\Big],\numberthis
\end{align*}
with $\eta_{r,\bm{v}} \coloneq \eta - \delta_{(r, \bm{v})}$. Notice we do not have $\m_{n}(\eta)\oplus r\bm{v} = \m_{n}(\eta + \delta_{(r, \bm{v})})$ in general, because $(r, \bm{v})$ may not be one of the $n$ first points of $\eta$; in that case $\m_{n}(\eta)\oplus r\bm{v} = \m_{n}(\eta)$. To deal with this difficulty, we will need the following lemma.

\begin{lemma}
Define
\begin{align*}
&\eta_{<} \coloneq \Big\lb (r, \bm{v}) \in \eta,\ r\bm{v} \in \big[\bm{0}, \m_{n}(\eta)\big)\Big\rb\\
&\eta_{=} \coloneq \Big\lb (r, \bm{v}) \in \eta,\ \exists j \in \lbra 1,d \rbra,\ rv^{j} = \m^{j}_{n}(\eta) \neq 0\Big\rb\\
&\eta_{>} \coloneq \Big\lb (r, \bm{v}) \in \eta,\ r\bm{v} \notin \big[\bm{0}, \m_{n}(\eta)\big]\Big\rb.
\end{align*}
With the above notations, those three sets are a.s. disjoint and partition $\eta$:
\begin{align}\label{3_dcpt_eta}
\eta = \eta_{<}\uplus \eta_{=} \uplus \eta_{>}.
\end{align}
\end{lemma}

\begin{proof}
Equality \eqref{3_dcpt_eta} indeed constitutes a partition of $\eta$ a.s.: if a vector $r\bm{v}$ does not have all of its coordinates strictly less than the corresponding coordinates of $\m_{n}(\eta)$, then it means that one of it coordinates is equal or greater than its counterpart of $\m_{n}(\eta)$. Those two cases are incompatible since $\eta$ is a Poisson process with diffuse intensity measure on $\Epol$. Indeed, assume there exist $j_{1}, j_{2}$ such that
\[
rv^{j_{1}} = \m^{j_{1}}_{n}(\eta) \neq 0\ \text{and\ } rv^{j_{2}} > \m^{j_{2}}_{n}(\eta).
\]
The second assumption implies that $(r, \bm{v})$ does not belong to the $n$ first points of $\eta$. But then it would imply that one could find another $(r_{i}, \bm{v}_{i})$ in $\eta$, with $i \in \lbra 1,n \rbra$, such that 
\[ 
rv^{j_{1}} = r_{i}v_{i}^{j_{1}}.
\] 
Because the measure $\alpha r^{-(\alpha + 1)}\dif r$ is diffuse on $\R_{+}^{*}$, this situation can occur with positive probability only if both $v_{i}^{j_{1}}$ and $v^{j_{1}}$ are equal to $0$. This contradicts the fact that $rv^{j_{1}}$ is not null.
\end{proof}
%A way to sum up the contents of the previous lemma is to say that each coordinate of $\m_{n}(\eta)$ comes from exactly one $r\bm{v}$, which does not prevent it from containing other coordinates of $\m_{n}(\eta)$. 
With that partition, we decompose the sum in \eqref{3_LePage_3parts} into three pieces, denoted by $S_{<}, S_{=}$ and $S_{>}$ and analyze them.

1. The first part corresponds to $\eta_{<}$ and is the easiest to deal with:
\begin{align*}
\sum_{(r,\bm{v}) \in \eta_{<}} \big\langle r\bm{v}, \nabla g_{h}\big(\m_{n}(\eta_{r,\bm{v}})\oplus r\bm{v}\big)\big\rangle_{\m_{n}(\eta_{r,\bm{v}})} = 0,
\end{align*}
because each term in the inner product vanishes. Indeed, due to the indicator functions, the sum is null as soon as all coordinates of $r\bm{v}$ are less than those of $\m_{n}(\eta_{r,\bm{v}})$. But since $(r, \bm{v}) \in \eta_{<}$, we already know that no coordinate of $r\bm{v}$ contributes to $\m_{n}(\eta)$, and thus to $\m_{n}(\eta_{r,\bm{v}})$. Thus 
\[ 
\esp[S_{<}] = 0.
\]

2. Next we take care of the sum over $\eta_{>}$. First observe that if $(r, \bm{v}) \in \eta_{>}$, then it cannot belong to the $n$ first points of $\eta$ and is not taken into account by $\m_{n}$. So we have
\[
\m_{n}(\eta_{r,\bm{v}}) = \m_{n}(\eta).
\]
Besides, we have the inclusion
\[
\lb rv^{j} > \m^{j}_{n}(\eta) \rb = \bigcap_{k=1}^{n}\lb rv^{j} > r_{k}v^{j}_{k}\rb \subseteq \bigcap_{k=1}^{n}\lb v^{j} > v^{j}_{k}\rb.
\]
The inequality comes from noticing that since $r \in \eta_{>}$, so that it belongs to $(r_{i})_{i\geq n+1}$, we have $r < r_{k}$ if $k \in \lbra 1,n \rbra$, because the sequence $(r_{i})_{i \geq 1}$ has been sorted in decreasing order. Thus, $rv^{j} > \m_{n}^{j}(\eta)$ is possible only if $v^{j} > v^{j}_{k}$ for all $k \in \lbra 1,n \rbra$.
\begin{align*}
\vert \esp[S_{>}]\vert &= \Big\vert\esp\Big[\sum_{(r,\bm{v}) \in \eta_{>}} \big\langle r\bm{v}, \nabla g_{h}\big(\m_{n}(\eta_{r,\bm{v}})\oplus r\bm{v}\big)\big\rangle_{\m_{n}(\eta_{r,\bm{v}})}\Big]\Big\vert\\
		&= \Big\vert\esp\Big[\sum_{(r,\bm{v}) \in \eta_{>}} \big\langle r\bm{v}, \nabla g_{h}\big(\m_{n}(\eta)\oplus r\bm{v}\big)\big\rangle_{\m_{n}(\eta)}\Big]\Big\vert\\
		&\leq \sum_{j=1}^{d}\esp\Big[\sum_{(r,\bm{v}) \in \eta_{>}}rv^{j} \big\vert\partial_{j}g_{h}\big(\m_{n}(\eta)\oplus r\bm{v}\big)\big\vert\ind_{\lb rv^{j} > \m^{j}_{n}(\eta) \rb}\Big]\\
		&\leq \alpha\sum_{j=1}^{d}\esp\Big[\sum_{i = n+1}^{\infty}r_{i}\min\big(\alpha, (r_{i}v^{j}_{i})\big)^{\alpha})\ind_{\bigcap_{k=1}^{n}\lb v^{j}_{i} > v^{j}_{k}\rb}\Big]\\
		&\leq \alpha\sum_{j=1}^{d}\sum_{i = n+1}^{\infty}\esp\big[r_{i}\min(\alpha, r_{i}^{\alpha})\big]\prob\Big(\bigcap_{k=1}^{n}\lb v^{j}_{i} > v^{j}_{k}\rb\Big)\\
		&\leq d\alpha\frac{1}{n}\sum_{i = n+1}^{\infty}\esp[r_{i}^{\alpha+1}],
\end{align*}
thanks to inequality \eqref{2_g_h_inqlty}. The term $1/n$ comes from the basic bound:
\[
\prob(X_{1} > X_{2},\dots, X_{1} > X_{n}) \leq \frac{1}{n},
\]
where $X_{1},\dots, X_{n}$ are $n$ \iid random variables, because the marks $(\bm{v}_{i})_{i\geq 1}$ are \iid by definition. We have also used that the $v^{j}$ are always less than $\ind_{\lb v^{j} \neq 0 \rb}$, thanks to the assumption we made on the reference norm at the beginning on this section. Set $\Gamma_{i} \coloneq r_{i}^{-\alpha}$. By definition of the $r_{i}$, which are the points of a Poisson process with intensity $\alpha r^{-(\alpha + 1)}$, we know that $\Gamma_{i}$ has the gamma distribution $\Gamma(i, 1)$ with shape parameter $i$ and scale parameter $1$. Consequently:
\begin{align*}
\esp[r_{i}^{\alpha+1}] &= \esp\big[\Gamma^{-(1 + 1/\alpha)}_{i}\big]\\
		&= \frac{1}{\Gamma(i)}\int_{0}^{\infty}x^{-(1 + \frac{1}{\alpha})}x^{i-1}e^{-x}\dif x\\
		&= \frac{\Gamma(i - 1 - 1/\alpha)}{\Gamma(i)}\\
		&\leq e^{2}\Big(i - 1 - \frac{1}{\alpha}\Big)^{-(1 + \frac{1}{\alpha})}, 
\end{align*}
thanks to the inequality 
\begin{align}\label{3_Kevckic}
\frac{\Gamma(x)}{\Gamma(y)} \leq \frac{x^{x-1}}{y^{y-1}}e^{y - x},\ y>x>1
\end{align}
which is a special case of a bound proved in \cite{Kevckic71}. We have taken ${x = i-1-1/\alpha}$ and $y = i$, for $i \geq n+1$ and $n\geq 2$. This last bound is the general term of a convergent series. The integral test for convergence (\textit{i.e.} the Maclaurin-Cauchy test) tells us that:
\[
\sum_{i = n+1}^{\infty}\frac{1}{i^{1+\frac{1}{\alpha}}} \underset{n \to \infty}{\sim} \frac{1}{n^{\frac{1}{\alpha}}}\cdotp
\]
In particular, by bounding $n/(n-1-1/\alpha)$ by $2$, we find:
\[
\vert \esp[S_{>}] \vert = \Big\vert\esp\Big[\sum_{(r,\bm{v}) \in \eta_{>}} \big\langle r\bm{v}, \nabla g_{h}\big(\m_{n}(\eta_{r,\bm{v}})\oplus r\bm{v}\big)\big\rangle_{\m_{n}(\eta_{r,\bm{v}})}\Big]\Big\vert \leq 2de^{2}\alpha \frac{1}{n^{1 + \frac{1}{\alpha}}}\cdotp
\]

3. Finally, we treat the case of $\eta_{=}$. We forcefully make appear the second part of the generator times $\alpha$. A consequence of lemma \ref{3_dcpt_eta} is that it equals:
\[ 
\big\langle \m_{n}(\eta), \nabla g_{h}(\m_{n}(\eta))\big\rangle = \sum_{(r,\bm{v}) \in \eta_{=}} \big\langle r\bm{v}, \nabla g_{h}(\m_{n}(\eta))\big\rangle_{\m_{n}(\eta)}. 
\]
Therefore the error we commit by making this substitution is:
\begin{align*}
S_{=} - \big\langle \m_{n}(\eta), \nabla g_{h}&(\m_{n}(\eta))\big\rangle\\
		&= \sum_{(r,\bm{v}) \in \eta_{=}} r\Big[\big\langle \bm{v},\nabla g_{h}\big(\m_{n}(\eta_{(r,\bm{v})})\oplus r\bm{v}\big)\big\rangle_{\m_{n}(\eta_{r,\bm{v}})} - \big\langle \bm{v}, \nabla g_{h}(\m_{n}(\eta))\big\rangle_{\m_{n}(\eta)}\Big]\\
		&= \sum_{j=1}^{d}\sum_{(r,\bm{v}) \in \eta_{=}} rv^{j}\Big[\partial_{j}g_{h}\big(\m_{n}(\eta_{r,\bm{v}})\oplus r\bm{v}\big)\ind_{\lb rv^{j} > \m_{n}^{j}(\eta_{r, \bm{v}})\rb} - \partial_{k}g_{h}(\m_{n}(\eta))\ind_{\lb rv^{j} > \m_{n}^{j}(\eta)\rb}\Big]\\
		&= \begin{dcases}
			\sum_{j=1}^{d}\sum_{(r,\bm{v}) \in \eta_{=}} rv^{j}\Big[\partial_{j}g_{h}\big(\m_{n}(\eta_{r,\bm{v}})\oplus r\bm{v}\big)\ind_{\lb rv^{j} > \m_{n}^{j}(\eta_{r, \bm{v}})\rb} - \partial_{j}g_{h}(\m_{n}(\eta))\ind_{\lb rv^{j} > \m_{n}^{j}(\eta)\rb}\Big]\\
			0,
			\end{dcases}
\end{align*}
the first case occurring as soon as one coordinate of $r_{n+1}\bm{v}_{n+1}$ is greater than its counterpart of $\m_{n}(\eta)$, while both terms in the subtraction are equal if $r_{n+1}\bm{v}_{n+1} \in [\bm{0}, \m_{n}(\eta)]$. Thus the error is bounded by
\begin{align*}
\esp\Big[\sum_{(r,\bm{v}) \in \eta_{=}} &rv^{j}\Big\vert\partial_{j}g_{h}\big(\m_{n}(\eta_{r,\bm{v}})\oplus r\bm{v}\big)\ind_{\lb rv^{j} > \m_{n}^{j}(\eta_{r, \bm{v}})\rb} - \partial_{j}g_{h}(\m_{n}(\eta))\ind_{\lb rv^{j} > \m_{n}^{j}(\eta)\rb}\Big\vert\Big]\\
		&= \esp\Big[\sum_{(r,\bm{v}) \in \eta_{=}} rv^{j}\Big\vert\partial_{j}g_{h}\big(\m_{n}(\eta_{r,\bm{v}})\oplus r\bm{v}\big)\ind_{\lb rv^{j} > \m_{n}^{j}(\eta_{r, \bm{v}})\rb} - \partial_{j}g_{h}(\m_{n}(\eta))\ind_{\lb rv^{j} > \m_{n}^{j}(\eta)\rb}\Big\vert\ind_{A_{n}}\Big],
\end{align*}
where 
\[
A_{n} \coloneq \bigcup_{l=1}^{d}\bigcap_{k=1}^{n} \lb v^{l}_{n+1} > v^{l}_{k}\rb.
\]
Besides, we can bound $rv^{j}$ by $r_{1}$, and the partial derivatives by $\alpha$. Gathering those arguments, we find:
\begin{align*}
\big\vert\esp[S_{=}] - \esp\big[\big\langle \bm{Z}_{n}, \nabla g_{h}(\bm{Z}_{n})\big\rangle\big]\big\vert &\leq 2d\alpha\esp[r_{1}\ind_{A_{n}}]\\
		&\leq 2d^{2}\alpha\esp\big[r_{1}\ind_{\bigcap_{k=1}^{n}\lb v^{j} > v^{j}_{k}\rb}\big]\\
		&= 2d^{2}\alpha\esp[r_{1}]\prob\Big(\bigcap_{k=1}^{n}\lb v^{j} > v^{j}_{k}\rb\Big)\\
		&\leq 2d^{2}\alpha\Gamma\Big(1 - \frac{1}{\alpha}\Big)\frac{1}{n}
\end{align*}
where $r_{1}$ has the Fréchet distribution $\mathcal{F}(\alpha)$. The presence of the term $d^{2}$ comes from the double sum: the sum over $j$ has $d$ terms, and the sum over $\eta_{=}$ has a random number of terms which is bounded by $d$.
\end{proof}

Through a slightly finer analysis, it is possible to obtain a better rate of convergence in the previous result. This comes at the price of working with smoother test functions, namely functions in $\lipbis$

\begin{corollary}\label{3_rate_LePage_d2}
We use the notations of proposition \ref{3_rate_LePage_Wass}. 
There exists a constant $C_{\alpha}>0$ depending only on $d$ and $\alpha$ such that
\[ 
\mathrm{d}_{[2]}(\bm{Z}_{n}, \bm{Z}) \leq C_{\alpha}\frac{1}{n^{1+\frac{1}{\alpha}}},\ n\geq 2. 
\]
Furthermore we have that $C_{\alpha} = O\Big(\big(1-\frac{1}{\alpha}\big)^{-2}\Big)$ when $\alpha$ goes to $1^{+}$.
\end{corollary}

\begin{proof}
Because a doubly 1-Lipschitz function $h$ on $\R_{+}^{d}$ is $1$-Lipschitz and of class $\mathcal{C}^{1}$ by definition, all the arguments given in the proof of proposition \ref{3_rate_LePage_Wass} apply again. A careful examination shows that we lose the rate of convergence of $n^{-(1+1/\alpha)}$ when dealing with $S_{=}$. In other words we must bound more accurately
\begin{align}\label{3_to_bound}
\sum_{(r,\bm{v}) \in \eta_{=}} rv^{j}\big\vert\partial_{j}g_{h}\big(\m_{n}(\eta_{r,\bm{v}})\oplus r\bm{v}\big)\ind_{\lb rv^{j} > \m_{n}^{j}(\eta_{r, \bm{v}})\rb} - \partial_{j}g_{h}(\m_{n}(\eta))\ind_{\lb rv^{j} > \m_{n}^{j}(\eta)\rb}\big\vert\ind_{A_{n}}
\end{align}
We define two subsets of $\lbra 1,d \rbra$:
\begin{align*}
&J_{1} \coloneq \big\lb j \in \lbra 1,d \rbra,\ rv^{j} = \m^{j}_{n}(\eta) \big\rb\\
&J_{2} \coloneq \big\lb j \in \lbra 1,d \rbra,\ r_{n+1}v^{j}_{n+1} < rv^{j}\big\rb
\end{align*}
Recall that unless $v^{j} = 0$, we have $r_{n+1}v^{j}_{n+1} \neq rv^{j}$ almost surely. 
To make the rest of the proof clear, we distinguish all four possible cases, depending on whether $j$ belongs to $J_{1}$ and/or $J_{2}$, or not.

1. $j \in J_{1}\cap J_{2}$ - In that case, both indicator functions are equal to $1$. Now, because $\partial_{j}g_{h}$ is $C_{2, \alpha}$-Lipschitz, we have:
\begin{align*}
rv^{j}\big\vert\partial_{j}g_{h}\big(\m_{n}(\eta_{r,\bm{v}})\oplus r\bm{v}\big)\ind_{\lb rv^{j} > \m_{n}^{j}(\eta_{r, \bm{v}})\rb} - \partial_{j}g_{h}(\m_{n}(\eta))&\ind_{\lb rv^{j} > \m_{n}^{j}(\eta)\rb}\big\vert\\
		&\leq C_{2, \alpha}rv^{j}\Vert \m_{n}(\eta_{r,\bm{v}})\oplus r\bm{v} - \m_{n}(\eta)\Vert_{1}\\
		&\leq C_{2, \alpha}r_{1}\Vert \m_{n}(\eta_{r,\bm{v}})\oplus r\bm{v} - \m_{n}(\eta)\Vert_{1}\\
		&\leq 2dC_{2, \alpha}r_{1}r_{n+1}.
\end{align*}
The coordinates of the vector $\m_{n}(\eta_{r,\bm{v}})\oplus r\bm{v} - \m_{n}(\eta)$ are either null or a factor of $r_{n+1}$ by some $v^{j} \in [0,1]$, hence the last inequality. 
To compute the expectation of $r_{1}r_{n+1}$, we let $E_{1},\dots, E_{n+1}$ be $n+1$ \iid random variables with the exponential distribution $\mathcal{E}(1)$. Then 
\[
(r_{1}, r_{n+1}) \eqdis \big(E_{1}^{-\frac{1}{\alpha}}, (E_{1}+\dots+E_{n+1})^{-\frac{1}{\alpha}}\big) \leq \big(E_{1}^{-\frac{1}{\alpha}}, (E_{2} + \dots + E_{n+1})^{-\frac{1}{\alpha}}\big).
\]
Thus, by setting $\Gamma_{n} \coloneq \sum_{i=2}^{n+1}E_{i}$, we find
\begin{align*}
\esp[r_{1}r_{n+1}] &\leq \esp[r_{1}]\esp\Big[\Gamma_{n}^{-\frac{1}{\alpha}}\Big]\\
		&= \Gamma\Big(1 - \frac{1}{\alpha}\Big)\frac{\Gamma(n - \frac{1}{\alpha})}{\Gamma(n)}\\
		&\leq \Gamma\Big(1 - \frac{1}{\alpha}\Big)\Big(n - \frac{1}{\alpha}\Big)^{-\frac{1}{\alpha}}\\
		&\leq 2\Gamma\Big(1 - \frac{1}{\alpha}\Big)\frac{1}{n^{\frac{1}{\alpha}}},
\end{align*}
by using once again inequality \eqref{3_Kevckic}. This gives us the term $n^{-1/\alpha}$ while the indicator function $\ind_{A_{n}}$, which is independent of the $(r_{i})_{i\geq 1}$, yields the term $n^{-1}$ as before. 
%Notice summing this bound with respect to $j$ makes appear a factor $d^{2}$, while $C_{2, \alpha}$ also depends on $d$. Thus the final constant depends on $d^{4}$. 

2. $j \in J_{1}^{c}\cap J_{2}$ - If $j \notin J_{1}$, then the contribution of $rv^{j}$ to $\m_{n}^{j}(\eta)$ is ignored and thus the second indicator function in \eqref{3_to_bound} vanishes. So does the first indicator function; otherwise $rv^{j}$ would have to be greater than both $r_{n+1}v^{j}_{n+1}$ and $\m_{n}^{j}(\eta_{r, \bm{v}}) = \m^{j}_{n}(\eta)$. This contradicts $j \notin J_{1}$. Consequently, both indicator functions are null.

3. $j \in J_{1}\cap J_{2}^{c}$ - Under the assumption that $j \notin J_{2}$, the first indicator function is equal to $0$. It also implies that $rv^{j} < r_{n+1}v^{j}_{n+1} \leq r_{n+1}$. The second term is not null and bounded by a constant, so that
\begin{align*}
rv^{j}\big\vert\partial_{j}g_{h}\big(\m_{n}(\eta_{r,\bm{v}})\oplus r\bm{v}\big)\ind_{\lb rv^{j} > \m_{n}^{j}(\eta_{r, \bm{v}})\rb} - \partial_{j}g_{h}(\m_{n}(\eta))\ind_{\lb rv^{j} > \m_{n}^{j}(\eta)\rb}\big\vert \leq 2\alpha r_{n+1}
\end{align*}

4. $j \in J_{1}^{c}\cap J_{2}^{c}$ - As seen previously, $j \notin J_{1}$ is enough to make both indicator functions vanish.

To prove the estimate on the constant, recall from proposition that $C_{2, \alpha}$ depends on the integral of $\gamma_{t}^{-1/\alpha}$, and observe that:
\[ 
\int_{0}^{\infty}\gamma_{t}^{-\frac{1}{\alpha}}\dif t = \int_{0}^{1}\gamma_{t}^{-\frac{1}{\alpha}}\dif t + \int_{1}^{\infty}\gamma_{t}^{-\frac{1}{\alpha}}\dif t \leq \int_{0}^{1}t^{-\frac{1}{\alpha}}\dif t + \int_{1}^{\infty}\gamma_{t}^{-\frac{1}{\alpha}}\dif t = \Big(1-\frac{1}{\alpha}\Big)^{-1} + \int_{1}^{\infty}\gamma_{t}^{-\frac{1}{\alpha}}\dif t. 
\]
Because $\Gamma(x) \underset{x \to 0^{+}}{\sim} x^{-1}$, we see that the constant in the  bound for the case $j \in J_{1}\cap J_{2}$ is of order $(1-1/\alpha)^{-2}$, hence concluding the proof.
case 
\end{proof}

\begin{remark}
We make two observations: first the bound of theorem \ref{3_rate_LePage_d2} becomes better as $\alpha$ gets closer to $1$, but in exchange the constant $C_{\alpha}$ explodes. Second we had to resort to the distance $\mathrm{d}_{[2]}$ to obtain this rate. Unlike the Kolmogorov distance, it is not invariant by monotonous bijective transformations and so we cannot deduce rates of convergence when $\alpha \in (0,1]$. 
\end{remark}

We bring a partial solution to both problems by using proposition 2.4. proved in \cite{Gaunt23}:

\begin{corollary}
Let $\alpha\in \R_{+}^{*}$ and assume that the angular measure $\nu$ is such that $\bm{Z} \sim \mathcal{MS}(\alpha, \nu)$ has a bounded density with respect to the Lebesgue measure on $\R^{d}$. Then there exist a constant $C>0$ depending only $d$ such that the de Haan-LePage series $(\bm{Z}_{n})_{n\geq 1}$ satisfies:
\[ 
\mathrm{d}_{\mathrm{Kol}}(\bm{Z}_{n}, \bm{Z}) \leq C\Big(\frac{\log n}{n}\Big)^{\frac{2}{3}},\ n \geq 2. 
\]
\end{corollary}

\begin{proof}
Let $\bm{Z}' \sim \mathcal{MS}(\beta, \nu)$. Notice it has the same distribution as $\bm{Z}^{\alpha/\beta}$ for every positive $\beta$. With obvious notations, we denote by $\bm{Z}_{n}'$ the corresponding partial de Haan-LePage series. Its law is the same as $\bm{Z}_{n}^{\alpha/\beta}$.

Now, thanks to proposition 2.4. and corollary \ref{3_rate_LePage_d2}, we know that there exists a constant $C$ independent of $n$ and $\beta$ such that
\[
\mathrm{d}_{\mathrm{Kol}}(\bm{Z}_{n}', \bm{Z}') \leq C\Big(1 - \frac{1}{\beta}\Big)^{-\frac{2}{3}}\frac{1}{n^{\frac{1}{3}(1 + \frac{1}{\beta})}}
\]
for $n$ greater than $2$ and any $\beta \in (1, +\infty)$. The Kolmogorov distance is invariant under non-decreasing transformations and Because $\bm{x} \mapsto \bm{x}^{\alpha/\beta}$ is non-decreasing, so the left-hand side is also equal to $\mathrm{d}_{\mathrm{Kol}}(\bm{Z}_{n}, \bm{Z})$. In particular, it does not depend on $\beta$. Thus, taking $\beta^{-1} = 1 - (\log n)^{-1}$, we find the announced result.
\end{proof}

%\printbibliography

%%% Local Variables:
%%% mode: latex
%%% TeX-master: "fc_JFA"
%%% End:

\bibliographystyle{amsplain}
\bibliography{biblio}
\end{document}